\documentclass{article}
\usepackage[utf8]{inputenc}
\usepackage[T1]{fontenc}
\usepackage{amsmath}
\usepackage{amsfonts}
\usepackage{amsthm}
\usepackage{amssymb}
\usepackage{color}

\usepackage[margin=1in]{geometry}

\providecommand{\bysame}{\leavevmode\hbox to3em{\hrulefill}\thinspace}
\providecommand{\MR}{\relax\ifhmode\unskip\space\fi MR }

\providecommand{\href}[2]{#2}

\newcommand{\cL}{{\mathbb {L}}}

\newcommand{\bpf}{\begin{preuve}}
\newcommand{\epf}{ \end{preuve} \medskip}

\newcommand{\benum}{\begin{enumerate}}
\newcommand{\eenum}{\end{enumerate}}

\newcommand{\bitem}{\begin{itemize}}
\newcommand{\eitem}{\end{itemize}}

\newcommand{\brmq}{\begin{rmq}}
\newcommand{\ermq}{\end{rmq}}

\newcommand{\brmqs}{\begin{rmqs}}
\newcommand{\ermqs}{\end{rmqs}}

\newcommand{\bapp}{\begin{application}}
\newcommand{\eapp}{\end{application}}

\newcommand{\bapps}{\begin{applications}}
\newcommand{\eapps}{\end{applications}}

\newcommand{\bdefi}{\begin{definition}}
\newcommand{\edefi}{\end{definition}}

\newcommand{\beq}{\begin{equation}}
\newcommand{\eeq}{\end{equation}}

\def\bpm{\begin{pmatrix}}
\def\epm{\end{pmatrix}}

\newcommand{\bcas}{\begin{cases}}
\newcommand{\ecas}{\end{cases}}

\newcommand{\bex}{\begin{exemp}}
\newcommand{\eex}{\end{exemp}}

\newcommand{\bexs}{\begin{exemps}}
\newcommand{\eexs}{\end{exemps}}

\newcommand{\bprop}{\begin{proposition}}
\newcommand{\eprop}{\end{proposition}}

\newcommand{\bthm}{\begin{theoreme}}
\newcommand{\ethm}{\end{theoreme}}

\newcommand{\bcor}{\begin{corollaire}}
\newcommand{\ecor}{\end{corollaire}}

\newcommand{\blem}{\begin{lemme}}
\newcommand{\elem}{\end{lemme}}

\newcommand{\beqna}{\begin{eqnarray}}
\newcommand{\eeqna}{\end{eqnarray}}

\newcommand{\beqnas}{\begin{eqnarray*}}
\newcommand{\eeqnas}{\end{eqnarray*}}

\newcommand{\bs}{\backslash}

\def\Id{{\rm{Id}}} 
\def\id{{\rm{id}}} 
\def\st{{\rm{st}}}

\def\eps{\varepsilon}

\def\diag{{\rm diag}}
\def\supp{{\rm supp}}
\def\htt{{\rm ht}}
\def\vol{{\rm vol}}

\def\cB{{\mathcal B }}
\def\cC{{\mathcal C}}
\def\cD{{\mathcal D}}
\def\cE{{\mathcal E }}
\def\cF{{\mathcal  F}}
\def\cG{{\mathcal  G}}

\def\cK{{\mathcal  K}}
\def\cL{{\mathcal L }}

\def\cS{{\mathcal S }}
\def\cT{{\mathcal T}}

\def\cX{{\mathcal X}}

\newcommand{\bbN}{{\mathbb {N}}}
\newcommand{\bbQ}{{\mathbb {Q}}}
\newcommand{\bbR}{{\mathbb {R}}}
\newcommand{\bbS}{{\mathbb {S}}} 
\newcommand{\bbT}{{\mathbb {T}}} 

\newcommand{\bbZ}{\mathbb {Z}}

\def\un{{\mathbf{1}}}

\def\bfb{{\mathbf{b}}}

\def\bfe{{\mathbf{e}}}

\def\bfj{{\mathbf{j}}}

\def\bfk{{\mathbf{k}}}

\def\bfm{{\mathbf{m}}}

\def\bfp{{\mathbf{p}}}
\def\bfq{{\mathbf{q}}}
\def\bfQ{{\mathbf{Q}}}

\def\bfV{{\mathbf{V}}}
\def\bfv{{\mathbf{v}}}

\def\bfx{{\mathbf{x}}}

\def\bfy{{\mathbf{y}}}

\def\bfz{{\mathbf{z}}}

\def\bf\Sigma{{\mathbf{\Sigma}}}

\def\bdxi{{\boldsymbol{\xi}}}
\def\bdalp{{\boldsymbol{\alpha}}}

\newtheorem{theoreme}{Theorem}[section]
\newtheorem{lemme}[theoreme]{Lemma}
\newtheorem{definition}[theoreme]{Definition}
\newtheorem{proposition}[theoreme]{Proposition}
\newtheorem{corollaire}[theoreme]{Corollary}

\newenvironment{exemp}{\noindent{\bf Example. --- }}{\par}
\newenvironment{exemps}{\noindent{\bf Examples}\benum}{\eenum\par}
\newtheorem{rmq}[theoreme]{Remark}
\newtheorem{rmqs}[theoreme]{Remarks}

\newenvironment{preuve}{\noindent{\it Proof. --- }}
{\hfill\rule{1.3mm}{2mm}\par}
\newenvironment{application}{\noindent{\bf Application. --- }}{\par}
\newenvironment{applications}{\noindent{\bf Applications. ---
}\benum}{\eenum\par}

\theoremstyle{definition}

\title{Asymptotics of the free path length in the periodic Lorentz gas \\ in the Boltzmann-Grad limit
}

\author{Songzi Li\thanks{School of Mathematics, Renmin University of China, 59, Zhongguancun Da Jie, Beijing, 100872, China ({\sf sli@ruc.edu.cn}). Supported by NSFC No.~11901569 and by fund (No.2018030249) from Renmin University of China.}}
\date{\today}

\begin{document}
\maketitle

\begin{abstract}
We study the distribution of the free path length in the periodic Lorentz gas in the Boltzmann-Grad limit. For any unimodular Euclidean lattice and any smooth initial
velocity distribution, we establish an asymptotic formula for the distribution function of the free path length with a polynomial rate of convergence. 
The proof combines the dynamical interpretation of the free path length in terms of expanding translates of horospherical orbits on the space of affine lattices, due to Marklof--Str\"ombergsson, with the effective equidistribution results of Kleinbock--Margulis and Kim. 
This extends the asymptotic formula of Boca--Zaharescu for the planar case to arbitrary dimension.
\end{abstract}
\textbf{Keywords}: periodic Lorentz gas, Boltzmann-Grad limit\\
\textbf{Mathematics Subject Classification (2020)}: Primary 60K50, 60F05; Secondary 37D45.

\section{Introduction and Main results}

The Lorentz gas, first proposed by Lorentz \cite{Lor1905} in 1905 to model electron motion in the metal, characterizes the dynamics of a test particle moving in an infinite array of fixed obstacles. This model is now regarded as a fundamental model for studying deterministic diffusion processes. The Boltzmann-Grad limit of the Lorentz gas corresponds to the regime where the obstacle radii tend to zero. In the periodic setting, where the obstacles are positioned at points of a Euclidean lattice, Marklof-Str\"ombergsson established the Boltzmann-Grad limit of the periodic Lorentz gas and derived the asymptotic estimates for the limiting distribution \cite{M-Sannals1, M-Sannals2, M-Sgafa}.  

The conjecture on the Boltzmann-Grad limit of the Lorentz gas originated in Lorentz's paper~\cite{Lor1905}, predicting that the particle density in the Boltzmann-Grad limit should obey a linear Boltzmann equation. In the case of random scatterer configuration, Lorentz's conjecture was first confirmed by Gallavotti~\cite{Gal1969}, see also Spohn~\cite{Spohn78}, Boldrighini-Bunimovich-Sinai~\cite{BBS1983}. However, Golse~\cite{Gol06} demonstrated that the conjecture fails for the periodic configuration. 

The free path length, defined as the distance between two successive collisions, serves as a fundamental observable in the Lorentz gas. As the first step in the proof of the convergence of the periodic Lorentz gas, Marklof-Str\"ombergsson \cite{M-Sannals1} established the Boltzmann-Grad limit of the free path length in the sense of distribution. Prior to this, Caglioti-Golse~\cite{Cag-Golse03} and Boca-Zaharescu \cite{Boca-Zaha07} proved the Boltzmann-Grad limit of the free path length in the planar periodic case. Notably, Boca-Zaharescu \cite{Boca-Zaha07} further derived the asymptotic formula of the distribution function of the free path length.

In this paper, we extend Boca-Zaharescu's \cite{Boca-Zaha07} results to high dimension, establishing the asymptotic formula for the distribution function of the free path length. To introduce the result,  we start with the notations, which largely follow~\cite{M-Sannals1}. Let $X$ denote the space of unimodular Euclidean lattices ${\rm SL(d, \bbZ)} \bs {\rm SL(d, \bbR)}$. For a lattice $\cL = \bbZ^{d}M \in X$, $M \in {\rm SL(d, \bbR)}$, we place a spherical obstacle of radius $r$ at each lattice point. Let $\cB^{d}_{r}$ and $ \bbS^{d-1}_{1}$ denote the ball in $\bbR^{d}$ with radius $r$ and the unit sphere in $\bbR^{d}$, respectively. 

The test particle moves at the speed $\bfv$ with $\|\bfv\| =1$ in the configuration space $\cK_{r} =\bbR^{d} \setminus (\cL + \cB^{d}_{r})$. We assume the absence of external potentials, so that the particle travels along straight lines until elastic collisions with obstacles occur. The phase space is defined as $T^{1}(\cK_{r}) = \cK_{r} \times \bbS^{d-1}_{1}$, equipped with the Liouville measure $d\mu_{r}(\bfq, \bfv) = dvol_{\bbR^{d}}(\bfq)dvol_{\bbS^{d-1}_{1}}(\bfv)$, where $dvol_{\bbR^{d}}$ and $dvol_{\bbS^{d-1}_{1}}$ denote the volume measures on $\bbR^{d}$ and $\bbS^{d-1}_{1}$, respectively.

For an initial position $\bfx = (\bfq, \bfv) \in  T^{1}(\cK_{r})$, the free path length on the lattice $\cL$ is defined as 
\beqna
\tau_{r}(\bfx) = \tau_{r}(\bfq, \bfv) = \inf\{t > 0, \bfq + t\bfv \notin \cK_{r}\}. 
\eeqna
This quantity represents the first collision time of the particle starting from the position $\bfq$ with the velocity $\bfv$. In~\cite{M-Sannals1}, the authors proved the convergence of the free path length in the sense of distribution. 

\bthm[Marklof-Str\"ombergsson]\label{ms.1}
Fix a unimodular Euclidean lattice $\cL = \bbZ^{d}M$ in $\bbR^{d}$.  Let $\bfq \in \bbR^{d} \backslash \cL$ and $\lambda$ be a Borel probability measure on $\bbS^{d-1}_{1}$ which is absolutely continuous with respect to the Lebesgue measure. Then there exists a continuous probability density $\Phi$ on $\bbR_{>0}$, which is independent of $\cL$ and $\bfq$,  such that for every $\xi \geq 0$, 
\beqna\label{bg.limit.free.micro}
\lim_{r \rightarrow 0}\lambda(\bfv \in \bbS^{d-1}_{1}, r^{d-1}\tau_{r}(\bfq, \bfv) \geq \xi) = \int^{\infty}_{\xi}\Phi(x)dx.
\eeqna
\ethm
We refer to \cite{M-Sannals1} for explicit formulas of the limiting density $\Phi$. Before the work \cite{M-Sannals1}, Boca-Zaharescu~\cite{Boca-Zaha07} proved the Boltzmann-Grad limit of the free path length of planar periodic Lorentz gas. They considered the torus $\bbT^{2} =  \bbR^{2}/\bbZ^{2}$ and the configuration space given by $\bbT^{2}_{r} =  (\bbR^{2} \setminus (\bbZ^{2} + \cB^{2}_{r}))/\bbZ^{2}$. The phase space $T^{1}(\bbT^{2}_{r})= \bbT^{2}_{r} \times \cS^{1}_{1}$ is equipped with the normalized Lebesgue measure $\mu_{r}$. The following asymptotic formula of the distribution function was obtained.

\bthm[Boca-Zaharescu]\label{BZ}
For every $\xi > 0$ and $\delta > 0$, 
\beqna\label{est.D1}
 \mu_{r}(r\tau_{r} > \xi) = \int^{\infty}_{\xi}\Phi(x)dx + O(r^{\frac{1}{8} - \delta}),
\eeqna
as $r \rightarrow 0$, where $\Phi$ is the limiting density of the free path length.
\ethm
We mention that $\Phi$ in this result is the same as the limiting distribution in \eqref{bg.limit.free.micro} when $d=2$.

 We prove the following asymptotic formula for the distribution function of the free path length with respect to the radius $r$. 
\bthm\label{m.1}
Let $\cL  = \bbZ^{d}M$ be a unimodular Euclidean lattice, and $\bdalp = - \bfq M^{-1}$ for $\bfq \in \bbR^{d} \backslash \cL$. Let $\lambda$ be a probability measure on $\bbS^{d-1}_{1}$ with smooth density with respect to the Lebesgue measure. Then, for any $\xi > 0$, we have
\beqna\label{est.main.l}
\lambda(\bfv \in \bbS^{d-1}_{1}: r^{d-1}\tau_{r}(\bfq, \bfv) \geq \xi ) =  \int^{\infty}_{\xi}\Phi (x)dx + E(\bdalp, r)
\eeqna
as $r \rightarrow 0$, where $\Phi$ is the limiting distribution of the free path length as defined in \eqref{bg.limit.free.micro}. There exist positive constants $\bar{\delta}$, $\bar{\kappa}$ and $\bar{\lambda}$ depending on $d$,  such that $E(\bdalp, r)$ is given by
\beqnas
E(\bdalp, r) = \left\{
\begin{array}{ll}
O(\zeta(\bdalp, r^{-\frac{d-1}{2}})^{- \bar{\delta}}+ r^{\bar{\kappa}}), & \bdalp \in \bbR^{d}  \setminus \bbQ^{d}, \\
O(r^{\bar{\lambda}}), & \bdalp \in \bbQ^{d}.
\end{array}
\right.
\eeqnas
The function $\zeta: \bbR^{d} \times \bbR^{+} \rightarrow \bbN$  is defined in \eqref{def.zeta}. The implicit constants in $E(\bdalp, r)$ depend on $M$, $\bfq$, $\lambda$ and $d$. 
\ethm

Following  \cite{M-Sannals1}, we introduce the "macroscopic coordinates": $(\bfQ, \bfV) = (r^{d-1}\bfq, \bfv)$. The corresponding free path length is then defined by $\cT_{r}(\bfQ, \bfV) = r^{d-1}\tau_{r}(r^{-(d-1)}\bfQ, \bfV)$. In this macroscopic setting, the asymptotic formula below extends the result of Boca-Zaharescu to higher dimensions. We note that our result confirms the polynomial decay of the free path length distribution with respect to the radius, although the precise exponent $\beta$ is not made explicit.
\bcor\label{c.1}
Let $\cL$ be as in Theorem~\ref{m.1}, and let $\Lambda$  be a Borel probability measure on $T^{1}(\bbR^{d})$ with a smooth density $p(\bfx, \bfv)$ compactly supported in $\bfx$ with respect to the Lebesgue measure. Assume that 
\beqna\label{asp.c1}
\sup_{\bfv \in \bbS^{d-1}_{1}}\sum^{d}_{i=1}\|\partial^{\gamma}_{i}p(\cdot, \bfv)\|_{L^{1}(\bbR^{d})} < +\infty
\eeqna
for some $\gamma > d$. Then there exists a positive constant $\beta$ such that, for any $\xi > 0$, 
\beqna\label{est.main.II}
\Lambda((\bfQ, \bfV) \in T^{1}(r^{d-1}\cK_{r}): \cT_{r}(\bfQ, \bfV) \geq \xi ) =  \int^{\infty}_{\xi}\Phi (x)dx + O(r^{\beta})
\eeqna
as $r \rightarrow 0$, where $\Phi$ is as in Theorem~\ref{m.1}. 
\ecor

Our proof relies on the interplay between the distribution of the free path length of the periodic Lorentz gas and expanding translates of certain unipotent subgroups on the affine lattice space, which was first observed in~\cite{M-Sannals1}. This connection enables us to reduce statistical properties of the free path length to the ergodic behavior of the corresponding flows. The asymptotic formula (Theorem~\ref{m.1}) makes use of the effective equidistribution results of such flows. When $\bdalp$ is rational, the distribution of the free path length corresponds to the flow on $X = {\rm SL(d, \bbZ)} \bs {\rm SL(d, \bbR)}$ and $X_{q} = \Gamma_{q} \bs  {\rm SL(d, \bbR)}$, as a special case of $\Gamma \bs G$, where $G$ is a Lie group and $\Gamma$ is a lattice in $G$. The exponential mixing property for such systems was proved by Kleinbock-Margulis~\cite{KM96}. These dynamical systems have been extensively studied due to their connections to ergodicity, geometry and number theory, particularly Diophantine approximation. For further details on the dynamics of expanding translates of the unipotent subgroup, we refer to \cite{KM96, KM99, KM12} and references therein; for the applications in Diophantine geometry, see \cite{D1985, KM98, BG19}.  When $\bdalp$ is irrational, the distribution of free path length corresponds to expanding translates of certain orbits on the affine lattice space ${\rm SL(d, \bbZ)} \ltimes \bbZ^{d} \bs {\rm SL(d, \bbR)} \ltimes \bbR^{d}$. The effective equidistribution in this setting, a recent breakthrough by Kim~\cite{Kim1}, provides the necessary foundation for our work. 

The remainder of this paper is organized as follows. In Section 2,  we provide the details on the connections between the distribution of the free path length and the dynamics of expanding translates of certain unipotent subgroup on the affine lattice space. Section 3 extends the effective equidistribution to a more general family of functions. Section 4 is devoted to the proofs of Theorem~\ref{m.1} and Corollary~\ref{c.1}.

\section{Distribution of free path length}
We begin by introducing some preliminary facts, including the connection between the distribution of the free path length in the Boltzmann-Grad limit and the ergodicity theorem of certain flows on the (affine) lattice space. Such flows can be considered as a special case of expanding translates of certain orbits in the (affine) lattice space. We then present the effective equidistribution results, which serve as the key role in our proof.

\subsection{Convergence of free path length distribution}
We recall the ergodicity theorems proved in~\cite{M-Sannals1} and their connections with the distribution of the free path length in the Boltzmann-Grad limit. 
Let 
\beqnas
\hat{G} &=&  {\rm SL(d, \bbR)} \ltimes \bbR^{d},  \ \hat{\Gamma} = {\rm SL(d, \bbZ)} \ltimes \bbZ^{d}, \\ 
Y &=&  \hat{\Gamma} \bs \hat{G} . 
\eeqnas
On $\hat{G}$, we have the multiplication law $(M, \bdxi)(M', \bdxi') = (MM', \bdxi M' + \bdxi')$. We write $(M, \bdxi)M' = (M, \bdxi)(M', \mathbf{0})$ for brevity. Define the action  of $\hat{G}$ on $\bbR^{d}$ by $\bfy \mapsto \bfy(M, \bdxi) = \bfy M +\bdxi$. 

Let $\cL = \bbZ^{d}M$ be the unimodular lattice. Write the affine lattice $\cL_{\boldsymbol{\alpha}} = \cL - \bfq$ as
$$
\cL_{\boldsymbol{\alpha}} = (\bbZ^{d} + \boldsymbol{\alpha})M = \bbZ^{d}(M,  \boldsymbol{\alpha}M) = \bbZ^{d}(\un_{d},  \boldsymbol{\alpha})M,
$$
where $\bfq = -  \boldsymbol{\alpha} M$.  When $\boldsymbol{\alpha} = \frac{\bfp}{q}$, the affine lattice space can be identified as
\beqnas
X_{q} = \Gamma(q) \backslash {\rm SL(d, \bbR)},
\eeqnas
where $\Gamma(q) = \{\gamma \in {\rm SL(d, \bbZ)}, \gamma \equiv \un_{d} \ mod \ q\}$. In the special case where $\boldsymbol{\alpha} \in \bbZ^{d}$, we denote
$$
X = X_{1} = {\rm SL(d, \bbZ)} \bs {\rm SL(d, \bbR)}
$$ as the space of unimodular lattice. Let $\mu_{Y}$, $\mu_{q}$, $\mu_{X} = \mu_{1}$ be the normalized right-invariant Haar measures on $Y$, $X_{q}$ and $X$, respectively.

We present the following ergodicity theorem obtained in~\cite{M-Sannals1}.  The ergodicity theorem for the irrational case $\boldsymbol{\alpha} \in \bbR^{d} \setminus \bbQ^{d}$ (see Theorem 5.2, \cite{M-Sannals1}) is a special case of Shah's theorem~\cite{shah}, as a consequence of Ratner's theory. For the rational case $\boldsymbol{\alpha} = \frac{\bfp}{q}$, $\bfp \in \bbZ^{d}$, $q \in \bbZ_{+}$,  the ergodicity theorem (see Theorem 5.8, \cite{M-Sannals1}) is implied by the mixing property of a subgroup in ${\rm SL(d, \bbR)}$ on $\Gamma \bs {\rm SL(d, \bbR)}$, with $\Gamma = \Gamma(q)$ for $\boldsymbol{\alpha} = \frac{\bfp}{q}$ and $\Gamma = {\rm SL(d, \bbZ)}$ for $ \boldsymbol{\alpha} \in \bbZ^{d}$ (which corresponds to the case $q =1$), see \cite{EskMcM93}.

\bthm[Marklof-Str\"ombergsson]\label{thm.ergodicity.2}
Let $\lambda$ be a Borel probability measure on $\bbR^{d-1}$ which is absolutely continuous with respect to Lebesgue measure. Let
\beqna\label{def.n}
n_{-}(\bfx) = \left(
\begin{array}{ccc}
1 & \bfx \\
\mathbf{0^{t}} & \un_{d-1} \\
\end{array}
\right),
\eeqna
and 
\beqna\label{def.at}
A^{t} =\left(
\begin{array}{ccc}
e^{-t} & \mathbf{0} \\
\mathbf{0^{t}} & e^{\frac{t}{d-1}}\un_{d-1} \\
\end{array}
\right).
\eeqna

\bitem
{\item
For $\boldsymbol{\alpha}  \in \bbR^{d} \setminus \bbQ^{d}$, $M \in {\rm SL(d, \bbR)}$ and a bounded continuous function $f: Y \rightarrow \bbR$, one has
\beqna\label{ergodicity.1}
\lim_{t \rightarrow \infty} \int_{\bbR^{d-1}} f \big((\un_{d},  \bdalp)Mn_{-}(\bfx)A^{t} \big) d\lambda(\bfx) = \int_{Y} f(g) d\mu_{Y}(g).
\eeqna
}
{\item
For  $\boldsymbol{\alpha} = \frac{\bfp}{q}$ or $\boldsymbol{\alpha} \in \bbZ^{d}$ (for $q=1$), $M \in {\rm SL(d, \bbR)}$ and a bounded continuous function $f: X_{q} \rightarrow \bbR$, one has
\beqna\label{ergodicity.q}
\lim_{t \rightarrow \infty} \int_{\bbR^{d-1}} f \big(Mn_{-}(\bfx)A^{t} \big) d\lambda(\bfx) = \int_{X_{q}} f(g) d\mu_{q}(g).
\eeqna
}
\eitem
\ethm

A key observation in \cite{M-Sannals1} is that the limiting distribution of the free path length is implied by that of the counting function in the some lattice-counting problem, which can be turned into the dynamics of expanding translates of certain orbits on the affine lattice space. The connection between the distribution of the free path length and that of the corresponding flow on the affine lattice space is reformulated into the following proposition. 

Define
$$
\varTheta_{r}(\sigma) =  \{(\bfv, \bfy) \in \bbS^{d-1}_{1} \times \bbR^{d}, 0 \leq \|\bfy\| < \sigma^{d-1}, \bfy \in \bbR_{\geq 0}\bfv + r\cS^{d-1}_{\bfv, \perp}\},
$$
where $\cS^{d-1}_{\bfv, \perp} = \{\bfz_{\bfv, \perp}, \bfz \in \bbS^{d-1}_{1}, \bfz \cdot \bfv  > 0\}$, $\bfz_{\bfv, \perp} = \bfz - (\bfz \cdot \bfv) \bfv$. Let $E_{1}:  \bbS^{d-1}_{1} \rightarrow SO(d)$ be the matrix such that $\bfv E_{1}(\bfv) = \bfe_{1}$ for $\bfv \in \bbS^{d-1}_{1} \setminus \{-\bfe_{1}\}$, and $e^{t} = r^{-(d-1)}$ in \eqref{def.at}. 

Define 
\beqna\label{def.theta.t}
\nonumber  \Theta_{t}(\sigma) &=& \varTheta_{r}(\sigma) E_{1}(\bfv)A^{t}:= \{(\bfv, \bfy E_{1}(\bfv)A^{t}),  (\bfv, \bfy) \in  \varTheta_{r}(\sigma) \} \\
 &=&  \{(\bfv, \bfy) \in \bbS^{d-1}_{1} \times \bbR^{d}, 0 \leq \|\bfy\| < \sigma^{d-1}, \bfy \in \bbR_{ \geq 0}\bfe_{1} + e^{-\frac{t}{d-1}}\cS^{d-1}_{1, \perp}\}A^{t}
\eeqna
where $\cS^{d-1}_{1, \perp} = \{\bfz_{\perp}, \bfz \in \bbS^{d-1}_{1}, \bfz_{1} > 0\}$, $\bfz_{\perp} = \bfz - (\bfz \cdot \bfe_{1})\bfe_{1}$, and
\beqna\label{def.theta}
\Theta(\sigma) &=&  \{(\bfv, \bfy) \in \bbS^{d-1}_{1} \times \bbR^{d}, y_{1} \in [0, \sigma^{d-1}), \|(y_{2}, \dots, y_{d})\| \leq 1\}.
\eeqna

For irrational $\bdalp$, define
\beqna\label{def.ce.t}
\cE_{t}(\sigma) = \{(\bfv, g) \in \bbS_{1}^{d-1} \times Y, \bbZ^{d}  g \cap \Theta_{t}(e^{\frac{t}{d-1}}\sigma)|_{\bfv} = \emptyset \},
\eeqna
and
\beqna\label{def.ce}
\cE (\sigma) = \{(\bfv, g) \in \bbS_{1}^{d-1} \times Y, \bbZ^{d}g \cap \Theta(\sigma)|_{\bfv} = \emptyset \}.
\eeqna

For rational $\bdalp = \frac{\bfp}{q}$ (including  $\bdalp \in \bbZ^{d}$ when $q = 1$), define
\beqna\label{def.ce.t.q}
\cE_{t, q}(\sigma) = \{(\bfv, g) \in \bbS_{1}^{d-1} \times X_{q}, (\bbZ^{d} +  \frac{\bfp}{q})g \cap \Theta_{t}(e^{\frac{t}{d-1}}\sigma)|_{\bfv} = \emptyset \},
\eeqna
and
\beqna\label{def.ce.q}
\cE_{q}(\sigma) = \{(\bfv, g) \in \bbS_{1}^{d-1} \times X_{q}, (\bbZ^{d}+  \frac{\bfp}{q})g \cap \Theta(\sigma)|_{\bfv} = \emptyset \}.
\eeqna
We point out that our definitions of $\Theta_{t}(\sigma)$,  $\Theta(\sigma)$, $\cE_{t}(\sigma)$, $\cE (\sigma)$ and $\cE_{t, q}(\sigma)$, $\cE_{q}(\sigma)$ are equivalent to those in \cite{M-Sannals1} up to a linear transformation
\beqna\label{def.c}
C_{\sigma} = \left(
\begin{array}{ccc}
\sigma^{-(d-1)} & \mathbf{0} \\
\mathbf{0^{t}} & \sigma \un_{d-1} \\
\end{array}
\right).
\eeqna


\bprop\label{MS.equi}
Consider the lattice $\cL = \bbZ^{d}M$.  Let $\lambda$ be a Borel probability measure on $\bbS^{d-1}_{1}$, which is absolutely continuous with respect to the Lebesgue measure. Let $\bdalp = -  \bfq M^{-1}$.  When $\bdalp \in \bbR^{d} \setminus \bbQ^{d}$, we have
\beqna\label{2.2.3}
\nonumber \lim_{r \rightarrow 0}\lambda(\bfv, r^{d-1}\tau_{r}(\bfq, \bfv) \geq \xi) &=& \lim_{t \rightarrow \infty}  \int_{\bbS^{d-1}_{1}}\un_{\cE_{t}(\sigma)} (\bfv, (\un_{d}, \bdalp)ME_{1}(\bfv)A^{t}) d\lambda(\bfv)\\
& =& \int_{\bbS^{d-1}_{1} \times Y} \un_{\cE(\sigma)}(\bfv, g)d\lambda(\bfv)d\mu_{Y}(g),
\eeqna
where $\sigma = \xi^{\frac{1}{d-1}}$. When $\bdalp = \frac{\bfp}{q}$ or $\bdalp \in \bbZ^{d}$, we have
\beqna\label{2.2.4}
\nonumber \lim_{r \rightarrow 0}\lambda(\bfv, r^{d-1}\tau_{r}(\bfq, \bfv) \geq \xi) &=& \lim_{t \rightarrow \infty}  \int_{\bbS^{d-1}_{1}}\un_{\cE_{t, q}(\sigma)} (\bfv, ME_{1}(\bfv)A^{t}) d\lambda(\bfv) \\
& =& \int_{\bbS^{d-1}_{1} \times X_{q}} \un_{\cE_{q}(\sigma)}(\bfv, g)d\lambda(\bfv)d\mu_{q}(g).
\eeqna
\eprop

\bpf
We prove the case when $\bdalp \in \bbR^{d} \setminus \bbQ^{d}$. The rational case $\bdalp = \frac{\bfp}{q}$ and $\bdalp \in \bbZ^{d}$ can be proved in the same spirit.

By the observation (display $(4.2)$) in \cite{M-Sannals1}, one has
\beqna\label{key.1}
\nonumber& \int_{\bbS^{d-1}_{1}}\un_{\cE_{t}(\bar{\sigma}_{t})} (\bfv, (\un_{d}, \bdalp)ME_{1}(\bfv)A^{t}) d\lambda(\bfv) \\
\nonumber & \qquad \qquad \qquad  \leq  \lambda(\bfv, r^{d-1}\tau_{r}(\bfq, \bfv) \geq \xi) \\
 & \qquad \qquad \qquad  \qquad \qquad \qquad     \leq  \int_{\bbS^{d-1}_{1}}\un_{\cE_{t}(\tilde{\sigma}_{t})} (\bfv,  (\un_{d}, \bdalp)ME_{1}(\bfv)A^{t} ) d\lambda(\bfv),
\eeqna
where $\bar{\sigma}_{t} = \sigma (1 + \sigma^{-(d-1)}e^{-\frac{dt}{d-1}})^{\frac{1}{d-1}}$, $\tilde{\sigma}_{t} = \sigma (1 - \sigma^{-(d-1)}e^{-\frac{dt}{d-1}})^{\frac{1}{d-1}}$.
Notice that
\beqnas
 \lim_{t \rightarrow \infty}  \int_{\bbS^{d-1}_{1}}\un_{\cE_{t}(\sigma)} (\bfv, (\un_{d}, \bdalp)ME_{1}(\bfv)A^{t} ) d\lambda(\bfv) &=& \lim_{t \rightarrow \infty}  \int_{\bbS^{d-1}_{1}}\un_{\cE_{t}(\tilde{\sigma}_{t})} (\bfv, (\un_{d}, \bdalp)ME_{1}(\bfv)A^{t} ) d\lambda(\bfv)\\
&=&  \lim_{t \rightarrow \infty}  \int_{\bbS^{d-1}_{1}}\un_{\cE_{t}(\bar{\sigma}_{t})} (\bfv, (\un_{d}, \bdalp)ME_{1}(\bfv)A^{t}) d\lambda(\bfv).
\eeqnas
Thus \eqref{2.2.3} can be derived from the ergodicity theorem (Theorem 5.6 in \cite{M-Sannals1}, which is a corollary of the general version of Theorem~\ref{thm.ergodicity.2}), and the fact that the set $\lim_{t \rightarrow \infty} \overline{\sup \cE_{t}(\sigma_{t})} \setminus \lim_{t \rightarrow \infty}(\inf \cE_{t}(\sigma_{t}))^{\circ}$ has measure zero (Lemma 9.2, \cite{M-Sannals1}).

 \epf

\subsection{Effective equidistributions on the (affine) lattice spaces }

 In this section, we review the effective equidistribution result implied by Kleinbock-Margulis~\cite{KM96} for the spaces $X_{q} = \Gamma(q) \bs {\rm SL(d, \bbR)}$ and $X =   SL(d, \bbZ) \bs SL(d, \bbR)$, as well as the result by Kim~\cite{Kim1} on the affine lattice space $Y =  SL(d, \bbZ) \ltimes \bbZ^{d} \bs SL(d, \bbR) \ltimes \bbR^{d}$. Moreover, we adapt Kim's result to our settings. 

\subsubsection{Preliminaries: height and norms}
We consider a general setting as in \cite{KM96}. Fix $m, n \in \bbN$, set $d=m+n$ and define
\beqna\label{def.h}
H = \{u_{M} | u_{M} = \begin{pmatrix}
\Id_{m}  & M\\
0 & \Id_{n}
\end{pmatrix}, 
M \in M_{m, n}\},
\eeqna
which is expanding horospherical with respect to 
\beqna\label{def.aat}
a^{t} = \diag \{e^{-nt}\un_{m}, e^{mt}\un_{n}\}
\eeqna
for $t > 0$. Let $\mu_{H}$ be the Haar measure on $H$.
We set the height of $x \in X$ by
\beqnas
\htt(x) = \sup \{\|\bfv g\|^{-1}, x = \Gamma g, \bfv \in \bbZ^{d} \setminus \{0\}\},
\eeqnas
where $\|\cdot\|$ is the supremum norm of the vector. 
Let $K(R) = \{x \in X, \htt(x) \leq R\}$, then for all $R > 0$, $K(R)$ is compact due to Mahler's compact criterion. Moreover, 
\beqna\label{est.ht.2}
\mu_{X}(X \setminus K(R)) \ll R^{-d}.
\eeqna
For $f \in \cC_{c}^{\infty}(G)$, define the norm
\beqnas
\|f\|_{\cC^{k}} = \sum_{\deg(Z) \leq k}\|D_{Z}f\|_{L^{\infty}}.
\eeqnas
We recall the Sobolev norm introduced in \cite{Kim1}. For every $V \in \cG$, define the differential operator $D_{V}$ on $\cC^{\infty}_{c}(X_{q})$ by $D_{V}\phi(x) = \frac{d}{dt}|_{t=0}\phi(xe^{tV})$. Let $D_{Z}$ be the differential operator with respect to a monomial generated by the basis of the Lie algebra $\cG$. 
For $l \in \bbN$, $f \in \cC^{\infty}_{c}(X_{q})$, define the norm $\cS^{X_{q}}_{l}$ by
\beqna\label{def.norm.X}
\cS^{X_{q}}_{l}(f)^{2} &=& \sum_{\deg(Z) \leq l}\int_{X_{q}}|\htt(x)^{l}D_{Z}f(x)|^{2}d\mu_{q}(x).
\eeqna
Similarly, define the Sobolev norm $\cS^{Y}_{l}$ on $\cC^{\infty}_{c}(Y)$ for $l \in \bbN$, $f \in \cC^{\infty}_{c}(Y)$ by
\beqna\label{def.norm.Y}
\cS^{Y}_{l}(f)^{2} = \sum_{\deg(\hat{Z}) \leq l}\int_{Y}|\htt(\pi(y))^{l}\cD_{\hat{Z}}f(y)|^{2}d\mu_{Y}(y),
\eeqna
where $\pi: Y \rightarrow X$ is the natural projection, and $\hat{Z}$ is the monomial generated by the basis of the Lie algebra $\hat{\cG}$.  We mention the following properties of $\cS^{Y}$ due to \cite{EMV} and \cite{Kim1}: for $f \in \cC^{\infty}_{c}(Y)$, $l$ large enough and $\deg(\hat{Z}) \leq d+2$
\beqna\label{norm.1}
\|\cD_{\hat{Z}}f\|_{L^{\infty}(Y)} \ll \cS^{Y}_{l}(f).
\eeqna
For $g \in G$ and  $f \in \cC^{\infty}_{c}(Y)$, define $R_{g}f(y) = f(yg)$. Then for $l$ large enough we have
\beqna\label{norm.2}
\cS^{Y}_{l}(R_{g}f) \ll \|g\|^{\kappa l}\cS^{Y}_{l}(f),
\eeqna
where $\kappa$ is a constant and  the norm on $G$ is defined by $\|g\| = \max_{1 \leq i, j \leq d}\{|g_{ij}|, |(g^{-1})_{ij}|\}$. Moreover, 
\beqna\label{norm.3}
 \|f - R_{g}f\|_{L^{\infty}} \ll d_{G}(\id, g)\cS^{Y}_{l}(f),
 \eeqna
where $d_{G}$ is the distance function on $G$. 

Throughout this paper, we fix $k$ sufficiently large so that Theorem~\ref{thm.km} holds for  $\cS = \cS^{X}_{k}$ and Theorem~\ref{thm.kim} holds for $\cS = \cS^{Y}_{k}$, and the above properties are satisfied by $\cS^{Y}_{k}$.  

\subsubsection{Effective equidistributions}
The effective ergodicity theorem proved by Kleinbock-Margulis~\cite{KM96} implies the exponential mixing property of the expanding translates of the horospherical subgroup $H$ on $\Gamma(q) \bs SL(d, \bbR)$. See also Theorem 2.2, \cite{BG19} for the explicit dependence of the error term on the norms of the density functions.
\bthm[Kleinbock-Margulis]\label{thm.km}
Let $V \subset H$ be a fixed neighborhood of the identity in $H$ with smooth boundary and compact closure. Let $d\mu = p d\mu_{H}$ with $p \in \cC_{c}^{\infty}(H)$, $\supp \ p \subset V$ be a probability measure on $H$. For any compact set $L \subset X_{q}$ and $x_{0} \in L$, there exist $\lambda = \lambda(m, n) > 0$ and $T(L) \geq 0$ such that for any $f \in \cC^{\infty}_{c}(X_{q})$ and $t \geq T(L)$,
\beqna\label{effe}
\int_{H}f(x_{0}ha^{t})d\mu(h)  - \int_{X_{q}}f(x)d\mu_{q}(x) = O(e^{-\lambda t}\|p\|_{\cC^{k}}\cS(f)),
\eeqna
where the implied constant depends on $m$, $n$ and $V$.
\ethm

Kim~\cite{Kim1} extended Kleinbock-Margulis' result to the space of affine lattices, which implies an effective version of the irrational case in Theorem~\ref{thm.ergodicity.2}. Define a function $\zeta: \bbR^{d}   \setminus \bbQ^{d} \times \bbR^{+} \rightarrow \bbN$ given by 
\beqna\label{def.zeta}
\zeta(\bfb, T) = \min \{N \in \bbN, \min_{1 \leq |q| \leq N} |q\bfb|_{\bbZ} \leq \frac{N^{2}}{T}\},
\eeqna
where  $|\cdot|_{\bbZ}$ denotes the supremum distance from $\mathbf{0} \in \bbT^{d}$. Note that $\zeta(\bfb, \cdot)$ is non-decreasing, unbounded, and $\zeta(\bfb, T) \leq T^{\frac{d}{2d+1}}$. 

\bthm[Kim]\label{thm.kim}
Let $V \subset H$ be a fixed neighborhood of the identity in $H$ with smooth boundary and compact closure. There exists a constant $\delta = \delta(d) > 0$ such that
\beqna\label{effe.eqdis.aff.1}
\frac{1}{\mu_{H}(V)}\int_{V}f(y_{0}ha^{t}) d\mu_{H}(h) - \int_{Y}f d\mu_{Y} = O(\zeta(\bfb_{0}, e^{\frac{nt}{2}})^{-\delta}\cS(f)),
\eeqna
holds for $f \in \cC^{\infty}_{c}(Y)$,  $y_{0} = \hat{\Gamma}(\un_{d}, \bfb_{0})g_{0} \in Y$ with $g_{0} \in {\rm SL}(d, \bbR)$, $\bfb_{0} \in \bbT^{d}$ and $t \geq 0$ such that  $\|g_{0}\| \leq \zeta(\bfb_{0}, e^{\frac{nt}{2}})^{\delta}$. The implied constant depends on $V$, $d$. 
\ethm


We then derive the following result as a corollary of Theorem~\ref{thm.kim}. Let $B^{H}(id, r)$ be the ball centered at the identity in $H$ with radius $r$. 

\bcor\label{effe.aff.2}
Let $V \subset H$ be the same as in Theorem \ref{thm.kim} and assume that $B^{H}(id, R^{-1}) \subset V \subset B^{H}(id, R)$ for some $R \geq 1$.  Let $d\mu = p d\mu_{H}$ be a probability measure on $V$ with $p \in \cC^{\infty}(V)$, $\|p\|_{\cC^{k}} < e^{\frac{\kappa_{2}}{2}t}$ for some positive constant $\kappa_{2}$. There exist constants $\delta' =  \delta'(d)  > 0$, $c = c(m, n, k) > 0$ such that, 
\beqna\label{effe.eqdis.aff.2}
\int_{V}f(y_{0}ha^{t}) d\mu(h)  - \int_{Y}f d\mu_{Y}  = O(R^{c}\zeta(\bfb_{0}, e^{\frac{nt}{2}})^{-\delta'}(1 + \|p\|_{\cC^{k}})\cS(f)),
\eeqna
holds for $f \in \cC^{\infty}_{c}(Y)$, $y_{0} = \hat{\Gamma}(\un_{d}, \bfb_{0})g_{0} \in Y$ with $g_{0} \in {\rm SL}(d, \bbR)$, $\bfb_{0} \in \bbT^{d}$ such that  $\|g_{0}\| \leq \zeta(\bfb_{0}, e^{\frac{nt}{2}})^{\delta'}$. The implied constant depends on $d$, $m$, $n$.
\ecor

\bpf
The proof proceeds by adapting the argument of \cite{Kim1}. Let $\pi: Y \rightarrow X$ be the projection. For $f \in \cC^{\infty}_{c}(Y)$, define $\bar{f}(x) = \int_{\pi^{-1}(x)} f(y)dm_{\pi^{-1}(x)}(y)$ and $\tilde{f}(y) = f(y) - \bar{f}(\pi(y))$, where $m_{\pi^{-1}(x)}$ can be considered as the uniform measure on $\bbT^{d}$. We then have
\beqnas
\int_{X}\bar{f}d\mu_{X} = \int_{Y}f(y)d\mu_{Y}, \ \cS^{X}(\bar{f}) \leq \cS^{Y}(f),
\eeqnas
and
\beqnas
\int_{V}f(y_{0}ha^{t}) d\mu(h)= \int_{V}\bar{f}(\pi(y_{0}ha^{t})) d\mu(h) + \int_{V}\tilde{f}(y_{0}ha^{t}) d\mu(h). 
\eeqnas
By the same reasoning as in the proof of Corollary 2.4, \cite{BG19}, the exponential mixing \eqref{effe} also holds for $d\mu = pd\mu_{H}$, $p \in \cC^{\infty}(V)$ with some $\lambda' > 0$. By tracking the dependence on $R$ of the error term in \cite{BG19}, we derive that
\beqna\label{est.2.5.1}
\int_{V}\bar{f}(\pi(y_{0}ha^{t})) d\mu(h) = \int_{X}\bar{f}d\mu_{X} + O(R^{c}e^{-\lambda' t}\|p\|_{\cC^{k}}\cS(\bar{f})) = \int_{Y}fd\mu_{Y} + O(R^{c}e^{-\lambda' t}\|p\|_{\cC^{k}}\cS(f)), 
\eeqna
where $c = c(m, n, k) > 0$. 

We claim that the estimate for $\int_{V}\tilde{f}(y_{0}ha^{t}) d\mu(h) $ given in $(4.30)$ and $(4.39)$ of \cite{Kim1} remains valid. Namely, for any $t$ satisfying $\|p\|_{\cC^{k}} < e^{\frac{\kappa_{2}}{2}t}$, 
\beqna\label{est.2.5.2}
|\int_{V}\tilde{f}(y_{0}ha^{t}) d\mu(h)| \ll \rho^{\frac{d\kappa_{5}}{2k(d-1)}}\cS(f),
\eeqna
where 
$$
\rho = \max \{ e^{-\kappa_{4}t}, c^{-\frac{1}{2d}}_{11}\zeta(\bfb_{0}, e^{\frac{nt}{2}})^{-\frac{1}{2d}}\}.
$$
We refer to the original paper for the definitions of the constants $\kappa_{4}$, $\kappa_{5}$ and $c_{11}$. The proof of $(4.39)$ involves no dependence on $p$ or $R$, while $(4.30)$ uses only $(3.13)$ and the geometric bound $|I\setminus I_{int}| \ll \eps^{d}r^{-(d^{2}-1)}$. The Fourier decay $(4.36)$ of Proposition 4.10 in \cite{Kim1} holds verbatim for $d\mu = pd\mu_{H}$. This follows from the estimates in Section 3 and 4 of \cite{Kim1}, which can be verified to hold for $d\mu$ in place of $\frac{1}{\mu_{H}(V)}d\mu_{H}$. 


It remains to verify for $d\mu = pd\mu_{H}$, the weight comparability $(3.13)$, $(3.15)$ in \cite{Kim1}. For the density $\pi_{*}p$, the equidistribution~\eqref{effe} on $X$ holds with error $\|p\|_{\cC^{k}}\htt(x_{0})^{\kappa_{1}}e^{-\lambda t}$. We deal this term as $(2.7)$ in \cite{Kim1}: there exists some constant $0 \leq \kappa_{2} \leq \frac{1}{2}$, such that if $\htt(\pi(y_{0})) < e^{\frac{\kappa_{2}}{2}t}$, $\|p\|_{\cC^{k}} < e^{\frac{\kappa_{2}}{2}t}$ and $e^{-\kappa_{2}t} < r$, then $\pi_{*}\mu_{y,t}(B^{X}(x,r)) \asymp r^{d^{2}-1}$ for every $x \in K(C_{4}r^{\frac{1}{d}})$. The condition $\htt(\pi(y_{0})) < e^{\frac{\kappa_{2}}{2}t}$ is guaranteed by the initial condition on $g_{0}$. Hence $(3.13)$, $(3.15)$ hold for $d\mu = pd\mu_{H}$. 

Finally, combining \eqref{est.2.5.1} and \eqref{est.2.5.2} yields \eqref{effe.eqdis.aff.2} by choosing $\delta' > 0$ such that $e^{-\lambda' t} + \rho^{\frac{d\kappa_{5}}{2k(d-1)}} \ll \zeta(\bfb_{0}, e^{\frac{nt}{2}})^{-\delta'}$, which follows from $\zeta(\bdalp, e^{\frac{nt}{2}}) \leq e^{\frac{dnt}{2(2d+1)}}$ and $\rho \ll \zeta(\bfb_{0}, e^{\frac{nt}{2}})^{-\frac{1}{2d}}$ (the definition of $\rho$).
\epf

In this paper, we use the following effective equidistribution result with respect to $A^{t}$ by taking $y_{0} = (\un_{d}, \bdalp)M$, $n=d-1$, $m= 1$, $t_{1} = (d-1)t$, such that $A^{t_{1}} = a^{t}$.

\bcor\label{effe.aff.3}
Let $V$, $d\mu = pd\mu_{H}$ be the same as in Corollary \ref{effe.aff.2} with $\|p\|_{\cC^{k}} < e^{\frac{\kappa_{2}}{2(d-1)}t_{1}}$. There exist constants $\delta' =  \delta'(d)  > 0$, $c = c(m, n, k) > 0$ such that
\beqna\label{effe.eqdis.aff.3}
\int_{V}f((\un_{d}, \bdalp)MhA^{t_{1}}) d\mu(h)  - \int_{Y}f d\mu_{Y}  = O(R^{c}\zeta(\bdalp, e^{\frac{t_{1}}{2}})^{-\delta'}(1 + \|p\|_{\cC^{k}})\cS(f)),
\eeqna
holds for $f \in \cC^{\infty}_{c}(Y)$, $M \in {\rm SL}(d, \bbR)$, $\bdalp \in \bbR^{d} \setminus \bbQ^{d}$ such that $\|M\| \leq \zeta(\bdalp, e^{\frac{t_{1}}{2}})^{\delta'}$. The implied constant depends on $d$.
\ecor

\section{Effective equidistribution with respect to $E_{1}$}

%

In this section, we extend the effective equidistribution result to smooth, compactly supported functions on $\bbS_{1}^{d-1} \times Y$  (resp., $ \bbS_{1}^{d-1} \times X_{q}$). 
For a domain $D \subset \bbR^{d-1}$, define the norm on $\cC^{\infty}(D \times Y)$(resp., $\cC^{\infty}(D \times X_{q})$) as
\beqnas
\|f\|_{\cC^{1}(D, \cS)} = \sup_{\bfv \in D}\big( \cS(f(\bfv , \cdot)) + \cS(\partial_{\bfv}f(\bfv , \cdot)) \big).
\eeqnas

\bprop\label{prop2.10}
Consider the affine lattice $\cL_{\boldsymbol{\alpha}} =  \bbZ^{d}(\un_{d},  \boldsymbol{\alpha})M$. Let $\lambda$ be the same as in Theorem~\ref{m.1}. Let $f$ be a smooth, compactly supported function defined on $\bbS_{1}^{d-1} \times Y$ (resp., $\bbS_{1}^{d-1} \times X_{q}$) with $\|f\|_{\cC^{1}(\bbS_{1}^{d-1}, \cS)} < \infty$, and $E_{1}: \bbS_{1}^{d-1} \setminus \{-\bfe_{1}\} \rightarrow SO(d)$ be the smooth map satisfying $\bfv E_{1}(\bfv) = \bfe_{1}$.
\bitem
{\item
For $\boldsymbol{\alpha}  \in \bbR^{d} \setminus \bbQ^{d}$, there exists a positive constant $\delta_{2}$ such that for $t$ large enough,
\beqna\label{est.2.10.1}
\nonumber &|\int_{\bbS^{d-1}_{1}} f(\bfv, (\un_{d}, \bdalp)ME_{1}(\bfv)A^{t}) d\lambda(\bfv) - \int_{\bbS^{d-1}_{1} \times Y} f(\bfv, g) d\lambda(\bfv)d\mu_{Y}(g)| \\
& \qquad \qquad \qquad \qquad \qquad  \leq C(d, \lambda)\|f\|_{\cC^{1}(\bbS^{d-1}_{1}, \cS)}(\zeta(\bdalp, e^{\frac{t}{2}})^{-\delta_{2}} + e^{- \frac{d}{2(d-1)} t}).
\eeqna
}
{\item
For  $\boldsymbol{\alpha} = \frac{\bfp}{q}$, there exists positive constant $\lambda_{2}$ such that for $t$ large enough,
\beqna\label{est.2.10.2}
|\int_{\bbS^{d-1}_{1}} f(\bfv, ME_{1}(\bfv)A^{t}) d\lambda(\bfv) - \int_{\bbS^{d-1}_{1} \times X_{q}} f(\bfv, g)d\lambda(\bfv) d\mu_{q}(g)|  \leq C(d, \lambda)\|f\|_{\cC^{1}(\bbS^{d-1}_{1}, \cS)}e^{-  \lambda_{2} t}.
\eeqna
}
\eitem
\eprop

%

We start with the following estimates. 
\blem\label{l32}
Let $\tilde{\lambda}$ be a probability measure on a compact subset $\tilde{V} \subset \bbR^{d-1}$ with smooth density. Let $f$, $\{f_{t}\}$ be smooth, compactly supported functions defined on $\bbR^{d-1} \times Y$ (resp., $\bbR^{d-1} \times X_{q}$) such that $\lim_{t \rightarrow \infty}f_{t} = f$ uniformly. Assume that there exists a nonnegative function $\epsilon(t): (0, \infty) \rightarrow \bbR^{+}$ satisfying $\lim_{t \rightarrow \infty} \epsilon(t) = 0$ and
\beqna\label{est.t}
\|f_{t} - f\|_{\cC(\bbR^{d-1} \times Y)(resp., \cC(\bbR^{d-1} \times X_{q}))} \leq \epsilon(t).
\eeqna
\bitem
{\item
For $\boldsymbol{\alpha}  \in \bbR^{d} \setminus \bbQ^{d}$, there exists $\delta_{1} > 0$ such that for $t$ large enough, 
\beqna\label{est.3.2.1}
& \left| \int_{\tilde{V}} f_{t}(\bfx, (\un_{d}, \bdalp)Mn_{-}(\bfx)A^{t}) d\tilde{\lambda}(\bfx) - \int_{\tilde{V} \times Y} f(\bfx, g) d\tilde{\lambda}(\bfx)d\mu_{Y}(g) \right|  \nonumber  \\
& \qquad  \qquad  \qquad   \ll \big( (\vol(\tilde{V})\|\frac{d\tilde{\lambda}}{d\bfx}\|_{\cC^{k}} + 1)\|f\|_{\cC(\tilde{V},  \cS)} + \|f\|_{\cC^{1}(\tilde{V}, \cC)})\big)\zeta(\bdalp, e^{\frac{t}{2}})^{-\delta_{1}} +  \epsilon(t).
\eeqna
}
{\item
For  $\boldsymbol{\alpha} = \frac{\bfp}{q}$, there exists $\lambda_{1} > 0$ such that for $t$ large enough
\beqna\label{est.3.2.2}
& \left| \int_{\tilde{V}} f_{t}(\bfx, Mn_{-}(\bfx)A^{t}) d\tilde{\lambda}(\bfx) - \int_{\tilde{V} \times X_{q}} f(\bfx, g)d\tilde{\lambda}(\bfx) d\mu_{q}(g) \right| \nonumber  \\
& \qquad  \qquad  \qquad   \ll \big( \vol(\tilde{V})\|\frac{d\tilde{\lambda}}{d\bfx}\|_{\cC^{k}}\|f\|_{\cC(\tilde{V},  \cS)} + \|f\|_{\cC^{1}(\tilde{V}, \cC)})e^{- \lambda_{1} t} +  \epsilon(t).
\eeqna
}
\eitem
\elem

\bpf
It suffices to prove \eqref{est.3.2.1}, since the estimate~\eqref{est.3.2.2} can be derived in the same spirit by making use of exponential mixing property (Theorem~\ref{thm.km}) on $ \Gamma_{q} \bs SL(d, \bbR)$.

Set $B_{r}(0) = B^{\bbR^{d-1}}(0, r)$. Let $\psi \in \cC^{\infty}_{c}(\bbR^{d-1})$ satisfying $0 \leq \psi \leq 1$, $\psi \equiv 1$ on $B_{\sqrt{d-1}}(0)$ and $\supp \ \psi \subset B_{2\sqrt{d-1}}(0)$. Define the set
$$
K = \{\bfk \in \bbZ^{d-1}, B_{4\sqrt{d-1}\delta}(\delta \bfk) \cap \tilde{V} \neq \emptyset\}. 
$$
Then $\tilde{V} \subset \bigcup_{\bfk \in K}B_{4\sqrt{d-1}\delta}(\delta \bfk)$, and 
\beqna\label{K}
\sharp K \asymp \vol(\tilde{V})\delta^{-(d-1)}.
\eeqna 
For $\bfx \in \bbR^{d-1}$, define 
$$
D(\bfx) = \sum_{\bfj \in \bbZ^{d-1}}\psi(\frac{1}{2\delta}(\bfx - \delta \bfj)), \qquad \psi_{\bfk}(\bfx) = D(\bfx)^{-1}\psi(\frac{1}{2\delta}(\bfx - \delta \bfk)).
$$
We then conclude that $1 \leq D(\bfx) \leq C_{0}(d)$, where $C_{0}(d) \geq 1$ is some positive constant depending on $d$, and $\supp \ \psi_{\bfk} \subset B_{4\sqrt{d-1}\delta}(\delta \bfk)$, 
\beqna\label{ps}
\sum_{\bfk \in \bbZ^{d-1}}\psi_{\bfk} \equiv 1, \qquad \|\psi_{\bfk}\|_{\cC^{k}} \ll \delta^{-k}.
\eeqna

Since $f$ is smooth and compactly supported, we have for all $\bfx \in B_{4\sqrt{d-1}\delta}(\bfx_{0}) \subset \bbR^{d-1}$, $\delta > 0$, 
\beqna
\label{uniform.1}|f(\bfx, g) - f(\bfx_{0}, g)| \leq 4\sqrt{d-1}\|\frac{\partial f}{\partial \bfx}(\cdot, g)\|_{\infty}\delta,
\eeqna
so that given \eqref{est.t}, 
\beqna\label{2.8.3}
|f_{t}(\bfx, g) - f(\bfx_{0}, g)| \ll \|f\|_{\cC^{1}(\tilde{V}, \cC)}\delta + \epsilon(t).
\eeqna
By \eqref{ps}, we have $\sum_{\bfk}\int_{\tilde{V}}\psi_{\bfk}(\bfx)d\tilde{\lambda}(\bfx) = 1$, and
\beqna\label{est.2.8.1}
\nonumber & &\int_{\tilde{V}} f_{t}(\bfx, (\un_{d}, \bdalp)Mn_{-}(\bfx)A^{t}) d\tilde{\lambda}(\bfx) =  \sum_{\bfk \in K}\int_{\tilde{V}} f_{t}(\bfx, (\un_{d}, \bdalp)Mn_{-}(\bfx)A^{t}) \psi_{\bfk}(\bfx) d\tilde{\lambda}(\bfx) \\
&=&  \sum_{\bfk \in K}\int_{\tilde{V}} f(\delta\bfk, (\un_{d}, \bdalp)Mn_{-}(\bfx)A^{t}) \psi_{\bfk}(\bfx) d\tilde{\lambda}(\bfx) + O(\|f\|_{\cC^{1}(\tilde{V}, \cC)}\delta +  \epsilon(t)).
\eeqna
Let $m_{\bfk} = \int_{\tilde{V}}\psi_{\bfk}(\bfx) d\tilde{\lambda}(\bfx)$. For $\bfx \in \supp \ \psi_{\bfk}$, we write $\bfy = \bfx - \delta\bfk$, then $\bfy \in B_{4\sqrt{d-1}\delta}(0)$. Define
\beqnas
\tilde{\psi}_{\bfk}(\bfy) &=& \psi_{\bfk}(\delta\bfk + \bfy), \\
g_{\bfk}(\bfy) &=& m^{-1}_{\bfk}\tilde{\psi}_{\bfk}(\bfy)\rho(\delta\bfk + \bfy),
\eeqnas
where $\rho = \frac{d\tilde{\lambda}}{d\bfx}$ is the density of $\tilde{\lambda}$. Then $\|\tilde{\psi}_{\bfk}(\bfy)\|_{\cC^{k}} \ll \delta^{-k}$, and we can verify that $g_{\bfk}$ is a probability density on $B_{4\sqrt{d-1}\delta}(0)$
$$
\int_{B_{4\sqrt{d-1}\delta}(0)}g_{\bfk}(\bfy)d\bfy = m^{-1}_{\bfk}\int_{B_{4\sqrt{d-1}\delta}(0)}\tilde{\psi}_{\bfk}(\bfy)\rho(\delta\bfk +  \bfy) d\bfy =  m^{-1}_{\bfk} \int_{\tilde{V}}\psi_{\bfk}(\bfx) d\tilde{\lambda}(\bfx) =1.
$$
Moreover,  $g_{\bfk} \in \cC_{c}^{\infty}(B_{4\sqrt{d-1}\delta}(0))$ and
 \beqna\label{gk}
 m_{\bfk}\|g_{\bfk}\|_{\cC^{k}} \ll \|\frac{d\tilde{\lambda}}{d\bfx}\|_{\cC^{k}} \delta^{-k}.
 \eeqna
Together with $n_{-}(\bfx) = n_{-}(\delta\bfk)n_{-}(\bfy)$, we obtain
\beqnas
\int_{\tilde{V}} f(\delta\bfk, (\un_{d}, \bdalp)Mn_{-}(\bfx)A^{t}) \psi_{\bfk}(\bfx) d\tilde{\lambda}(\bfx) =  m_{\bfk}\int_{B_{4\sqrt{d-1}\delta}(0)} f(\delta\bfk, (\un_{d}, \bdalp)Mn_{-}(\delta\bfk)n_{-}(\bfy)A^{t}) g_{\bfk}(\bfy) d\bfy.
\eeqnas
Notice that for $V = \{n_{-}(\bfy), \bfy \in B_{4\sqrt{d-1}\delta}(0)\}$ with $0 < \delta < \frac{1}{4\sqrt{d-1}}$, we have $B^{H}(\id, \delta) \subset V \subset B^{H}(\id, \delta^{-1})$. By \eqref{effe.eqdis.aff.3}, we derive that 
\beqna\label{2.8.2}
\nonumber & &|\int_{\tilde{V}} f(\delta\bfk, (\un_{d}, \bdalp)Mn_{-}(\bfx)A^{t}) \psi_{\bfk}(\bfx) d\tilde{\lambda}(\bfx)  - m_{\bfk}\int_{Y}f(\delta \bfk, g)d\mu_{Y}(g)|\\
\nonumber &=&m_{\bfk}|\int_{B_{4\sqrt{d-1}\delta}(0)} f(\delta\bfk, (\un_{d}, \bdalp)Mn_{-}(\delta\bfk)n_{-}(\bfy)A^{t}) g_{\bfk}(\bfy) d\bfy -  \int_{Y}f(\delta \bfk, g)d\mu_{Y}(g)|\\
&\leq& C(d)\delta^{-c}m_{\bfk}(\|g_{\bfk}\|_{\cC^{k}} + 1)\cS(f(\delta \bfk, \cdot))\zeta(\bdalp, e^{\frac{t}{2}})^{-\delta'}
\eeqna
holds for $\|Mn_{-}(\delta\bfk)\| \leq \zeta(\bdalp, e^{\frac{t}{2}})^{\delta'}$. Since $\tilde{V}$ is bounded, $\|Mn_{-}(\delta\bfk)\| \leq C(\tilde{V})\|M\|$ for $\bfk \in K$. Thus we can take $t$ large enough such that $\|Mn_{-}(\delta\bfk)\| \leq \zeta(\bdalp, e^{\frac{t}{2}})^{\delta'}$ holds for all $\bfk \in K$. Then by \eqref{K} and \eqref{gk}, we derive that
\beqna\label{est.2.8.2}
\nonumber & &|\sum_{\bfk \in K}\int_{\tilde{V}} f(\delta\bfk, (\un_{d}, \bdalp)Mn_{-}(\bfx)A^{t}) \psi_{\bfk}(\bfx) d\tilde{\lambda}(\bfx)  - \int_{\tilde{V} \times Y}f(\bfx, g)d\tilde{\lambda}(\bfx)d\mu_{Y}(g)|\\
\nonumber &\leq& C(d)\delta^{-c}( \sharp K \delta^{-k}\|\frac{d\tilde{\lambda}}{d\bfx}\|_{\cC^{k}}  + 1)\|f\|_{\cC(\tilde{V}, \cS)}\zeta(\bdalp, e^{\frac{t}{2}})^{-\delta'}  + \|f\|_{\cC^{1}(\tilde{V}, \cC)}\delta\\
&\leq& C(d)\delta^{-c}(\delta^{-k - (d-1)}\vol(\tilde{V})\|\frac{d\tilde{\lambda}}{d\bfx}\|_{\cC^{k}} + 1)\|f\|_{\cC(\tilde{V}, \cS)}\zeta(\bdalp, e^{\frac{t}{2}})^{-\delta'} + \|f\|_{\cC^{1}(\tilde{V}, \cC)}\delta. 
\eeqna
Inserting \eqref{est.2.8.2} into \eqref{est.2.8.1} yields
\beqnas
& &|\int_{\tilde{V}} f_{t}(\bfx, (\un_{d}, \bdalp)Mn_{-}(\bfx)A^{t}) d\tilde{\lambda}(\bfx) - \int_{\tilde{V} \times Y}f(\bfx, g)d\tilde{\lambda}(\bfx)d\mu_{Y}(g)| \\
&\ll& \delta^{-c - k - d+1}(\vol(\tilde{V})\|\frac{d\tilde{\lambda}}{d\bfx}\|_{\cC^{k}} + 1)\|f\|_{\cC(\tilde{V}, \cS)}\zeta(\bdalp, e^{\frac{t}{2}})^{-\delta'} + \|f\|_{\cC^{1}(\tilde{V}, \cC)}\delta+  \epsilon(t).
\eeqnas 
Thus we prove \eqref{est.3.2.1} by taking $\delta^{c + d + k} = \zeta(\bdalp, e^{\frac{t}{2}})^{-\delta'} $ and $\delta_{1} = \frac{1}{c + d +k}\delta'$. 
\epf

We now utilize the observations in \cite{M-Sannals1} to prove Proposition~\ref{prop2.10}.

\bpf
We first extend the estimates with respect to $n_{-}(\bfx)$ in Lemma~\ref{l32} to that with respect to $E_{1}(\bfv)$. As pointed out in \cite{M-Sannals1}, one can define $E_{1}(\bfv)$ as following (which is equivalent to the one mentioned in \cite{M-Sannals1}),
\beqnas
E_{1}(\bfv) = \left(
\begin{array}{cc}
v_{1} & -\bfv_{\perp} \\
\bfv_{\perp}^{tr} & \Id_{d-1} - \frac{\bfv_{\perp}^{tr}\bfv_{\perp}}{1 + v_{1}}
\end{array}\right)
\eeqnas
such that it is smooth on $\bbS^{d-1}_{1} \setminus \{- \bfe_{1}\}$ and $E_{1}(\bfe_{1}) = \Id$. An observation in \cite{M-Sannals1} is that
\beqnas
E_{1}(\bfv)A^{t} =  n_{-}(\tilde{\bfv})A^{t}E_{2}(\bfv),
\eeqnas
where $\tilde{\bfv} = F(\bfv) = -v_{1}^{-1}\bfv_{\perp}$ and
\beqna\label{def.e3}
E_{2}(\bfv) =
\left(
\begin{array}{cc}
v_{1}^{-1} & 0 \\
\bfv_{\perp}^{tr}e^{- \frac{d}{d-1}t} & \Id_{d-1} - \frac{\bfv_{\perp}^{tr}\bfv_{\perp}}{1 + v_{1}}
\end{array}\right).
\eeqna
Notice that  $E_{2}(\bfv) \in SL(d, \bbR)$, $E_{2}(\bfv) \notin SO(d)$. Since $F$ has nonsingular differential at $\bfe_{1}$, there is an open set $D_{0}$ around $\bfe_{1}$, such that $F: D_{0} \rightarrow \tilde{D}_{0} \subset \bbR^{d-1}$ is a diffeomorphism. Let $B \subset D_{0}$ be a Borel set with $\lambda(B) > 0$ and $\tilde{B} \subset \tilde{D}_{0}$ be its image under $F$. Let $\tilde{\lambda}$ be the push-forward of the normalized measure $\lambda(B)^{-1}\lambda|_{B}$, i.e., $\tilde{\lambda}(\tilde{\bfv}) := \lambda(B)^{-1}\lambda|_{B} \cdot F^{-1}(\tilde{\bfv})$, then $\tilde{\lambda}$ is a probability measure on $\tilde{B}$. 
%
%
Since $E_{2}(\bfv)$ depends on $t$, we define $\tilde{f}_{t}, \tilde{f}: \bbR^{d-1} \times Y \rightarrow \bbR$ by
\beqna\label{def.til.ft}
\tilde{f}_{t}(\tilde{\bfv}, g) =  f(\bfv, gE_{2}(\bfv)),
\eeqna
and
\beqna\label{def.til.f}
\tilde{f}(\tilde{\bfv}, g) = f(\bfv, gE_{3}(\bfv)), \quad E_{3}(\bfv) = \left(
\begin{array}{cc}
v_{1}^{-1} & 0 \\
0 &  \Id_{d-1} - \frac{\bfv_{\perp}^{tr}\bfv_{\perp}}{1 + v_{1}}
\end{array}\right).
\eeqna
Then $\tilde{f}_{t}$, $\tilde{f}$ are smooth and compactly supported. It follows that the convergence $\lim_{t \rightarrow \infty}\tilde{f}_{t}(\tilde{\bfv}, g) = \tilde{f}(\tilde{\bfv}, g)$ holds uniformly on compact sets in $\bbR^{d-1} \times Y$. Notice that
\beqnas
E_{2}(\bfv) = E_{3}(\bfv)E_{4}(\bfv), \quad E_{4}(\bfv)=
\left(\begin{array}{cc}
1 & 0 \\
A(\bfv)e^{- \frac{d}{d-1}t} & \Id_{d-1}
\end{array}\right),
\eeqnas
where $A(\bfv) = (\Id_{d-1} - \frac{\bfv_{\perp}^{tr}\bfv_{\perp}}{1 + v_{1}})^{-1}\bfv_{\perp}^{tr} = (\Id_{d-1} + \frac{\bfv_{\perp}^{tr}\bfv_{\perp}}{v_{1}(1 + v_{1})})\bfv_{\perp}^{tr}$. 
Thus by \eqref{norm.3} we have for $\tilde{\bfv} \in \tilde{B}$
\beqna\label{2.8.4}
\nonumber |\tilde{f}_{t}(\tilde{\bfv}, g) - \tilde{f}(\tilde{\bfv}, g) | &\ll& \|f\|_{\cC(B, \cS)}d_{G}(E_{4}(\bfv), \Id_{d})\\
&\ll& \|f\|_{\cC(B, \cS)}\|A(\bfv)\|e^{- \frac{d}{d-1}t}.
\eeqna
Due to \eqref{norm.2},  we have
\beqnas
\cS(\tilde{f}(\bfv, \cdot)) &\ll& \|E_{3}(\bfv)\|^{\kappa k}\cS(f(\bfv, \cdot)), \\ 
  \|D_{\tilde{\bfv}}\tilde{f}\|_{\infty} &\ll& \|D_{\bfv}f \|_{\infty} + \|D_{\bfv} E_{3}\|_{\infty}\cS(f(\bfv, \cdot)).
 \eeqnas 
 Moreover, 
\beqna\label{2.8.5}
\int_{\tilde{B}}\tilde{f}_{t}(\tilde{\bfv}, (\un_{d}, \bdalp) Mn_{-}(\tilde{\bfv})A^{t})d\tilde{\lambda}(\tilde{\bfv}) =\lambda(B)^{-1} \int_{B}f(\bfv, (\un_{d}, \bdalp)ME_{1}(\bfv)A^{t})d\lambda(\bfv), 
\eeqna
and
\beqna\label{2.8.6}
\nonumber \int_{\tilde{B}  \times Y}\tilde{f}(\tilde{\bfv}, g)d\tilde{\lambda}(\tilde{\bfv})d\mu_{Y}(g) &=&\lambda(B)^{-1} \int_{B \times Y}f(\bfv, gE_{3}(\bfv))d\lambda(\bfv)d\mu_{Y}(g),\\
&=&\lambda(B)^{-1} \int_{B \times Y}f(\bfv, g)d\lambda(\bfv)d\mu_{Y}(g),
\eeqna
where the last line is due to the invariance of the Haar measure $\mu_{Y}$. For $M: \bbS^{d-1}_{1} \rightarrow G$, let $\|M\|_{B} =  \sup_{\bfv \in B} \|M(\bfv)\|$ for $B \subset \bbS^{d-1}_{1}$.  Applying \eqref{est.3.2.1} to $\tilde{f}_{t}$, $\tilde{f}$ with \eqref{2.8.4}, we obtain by \eqref{2.8.5} and \eqref{2.8.6} that
\beqna\label{est.4.3.0}
\nonumber & &|\int_{B}f(\bfv, (\un_{d}, \bdalp)ME_{1}(\bfv)A^{t})d\lambda(\bfv) - \int_{B \times Y}f(\bfv, g)d\lambda(\bfv)d\mu_{Y}(g)|\\
\nonumber &=&\lambda(B) |\int_{\tilde{B}}\tilde{f}_{t}(\tilde{\bfv}, (\un_{d}, \bdalp)Mn_{-}(\tilde{\bfv})A^{t})d\tilde{\lambda}(\tilde{\bfv}) - \int_{\tilde{B} \times Y}\tilde{f}(\tilde{\bfv}, g)d\tilde{\lambda}(\tilde{\bfv})d\mu_{Y}(g)| \\
\nonumber &\ll& \big((\vol(\tilde{B})\|\frac{d\tilde{\lambda}}{d\bfx}\|_{\cC^{k}} + 1)\|E_{3}\|_{B}^{\kappa k} \|f\|_{\cC(B,  \cS)} + (1 + \|D_{\bfv}E_{3}\|_{B})\|f\|_{\cC^{1}(B, \cS)} \big) \zeta(\bdalp, e^{\frac{t}{2}})^{-\delta_{1}} \\
& &+ \|A\|_{B}\|f\|_{\cC(B, \cS)}e^{- \frac{d}{d-1}t}.
\eeqna 



We now extend the estimates from $B$ to $\bbS^{d-1}_{1}$. Notice that $F$ is unbounded when $v_{1} \rightarrow 0$. Let $K_{r} = \bbS^{d-1}_{1} \cap (\cB^{d}_{-\bfe_{1}}(r) \cup  \{|v_{1}| \leq r\})$ for $0 < r < \frac{1}{2}$. Then $\lambda(K_{r}) \ll r$. For $\bfv \in \bbS^{d-1}_{1}\setminus K_{r}$, $\|F(\bfv)\| \leq r^{-1}$. Moreover, for $X \in \bbR^{d-1}$
\beqnas
F^{-1}(X) = (\frac{1}{\sqrt{1 + \|X\|^{2}}}, -\frac{X}{\sqrt{1 + \|X\|^{2}}}).
\eeqnas
or $F^{-1}(X) = (-\frac{1}{\sqrt{1 + \|X\|^{2}}}, \frac{X}{\sqrt{1 + \|X\|^{2}}})$. 
One can check that $\|F^{-1}\|_{\cC^{k}}$ is uniformly bounded on $F(K^{c}_{r}) \subset \bbR^{d-1}$. However, $\tilde{\lambda}$ is not smooth when restricted on $F(K^{c}_{r})$. To this end, let $\{h_{\eps}\}_{\eps < r}$ be the smooth approximation of $\un_{K^{c}_{r}}$, such that $\un_{K^{c}_{r}} \leq h_{\eps} \leq 1$, $\supp \ h_{\eps} \subset K^{c}_{r, \eps} := \{\bfv \in \bbS^{d-1}_{1}, d(\bfv, K^{c}_{r}) < \eps\}$ and $\|h_{\eps}\|_{\cC^{k}} \ll \eps^{-k}$. Let $\lambda_{\eps} =  \lambda h_{\eps}$. We introduce the mollified version of $\tilde{\lambda}_{\eps}$, given by 
\beqnas
\tilde{\lambda}_{\eps} = \frac{1}{Z_{r, \eps}} \lambda_{\eps} \cdot F^{-1}, 
\eeqnas 
where $Z_{r, \eps} = \lambda_{\eps}(K^{c}_{r, \eps})$. Then, $\|\frac{d\tilde{\lambda}_{\eps}}{d\bfx}\|_{\cC^{k}} \ll \eps^{-k}\|\frac{d\lambda}{d\bfx}\|_{\cC^{k}}$, $\vol(F(K^{c}_{r, \eps})) \asymp (r -\eps)^{-(d-1)}$ and 
\beqnas
\|E_{3}\|_{K^{c}_{r, \eps}} \ll (r - \eps)^{-1}, \ \|D_{\bfv}E_{3}\|_{K^{c}_{r, \eps}} \ll (r - \eps)^{-3}, \ \|A\|_{K^{c}_{r, \eps}} \ll (r - \eps)^{-1}, 
\eeqnas 
such that
\beqnas
& &|\int_{K^{c}_{r}}f(\bfv, (\un_{d}, \bdalp)ME_{1}(\bfv)A^{t})d\lambda(\bfv) - \int_{K^{c}_{r} \times Y}f(\bfv, g)d\lambda(\bfv)d\mu_{Y}(g)| \\
&\leq& |\int_{K^{c}_{r, \eps}}f(\bfv, (\un_{d}, \bdalp)ME_{1}(\bfv)A^{t})d\lambda_{\eps}(\bfv) - \int_{K^{c}_{r, \eps} \times Y}f(\bfv, g)d\lambda_{\eps}(\bfv)d\mu_{Y}(g)| + 2\|f\|_{L^{\infty}}\eps\\
&\ll& \big((1 + (r -\eps)^{-(d-1)}\eps^{-k}\|\frac{d\lambda}{d\bfx}\|_{\cC^{k}})(r - \eps)^{-\kappa k} + 1 + (r - \eps)^{-3})\|f\|_{\cC^{1}(\bbS^{d-1}_{1}, \cS)}\zeta(\bdalp, e^{\frac{t}{2}})^{-\delta_{1}} \\
& &+ (r - \eps)^{-1}\|f\|_{\cC(\bbS^{d-1}_{1}, \cS)}e^{- \frac{d}{d-1}t}+ \|f\|_{L^{\infty}}\eps,
\eeqnas
where the last line is due to \eqref{est.4.3.0}. Taking $\eps = \frac{r}{2}$, we conclude that
\beqna\label{est.4.3.1}
\nonumber & &|\int_{\bbS^{d-1}_{1}}f(\bfv, (\un_{d}, \bdalp)ME_{1}(\bfv)A^{t})d\lambda(\bfv) - \int_{\bbS^{d-1}_{1} \times Y}f(\bfv, g)d\lambda(\bfv)d\mu_{Y}(g)|\\
\nonumber &\leq&|\int_{K^{c}_{r}}f(\bfv, (\un_{d}, \bdalp)ME_{1}(\bfv)A^{t})d\lambda(\bfv) - \int_{K^{c}_{r} \times Y}f(\bfv, g)d\lambda(\bfv)d\mu_{Y}(g)| + |\int_{K_{r}}f(\bfv, (\un_{d}, \bdalp)ME_{1}(\bfv)A^{t})d\lambda(\bfv) |\\
\nonumber & & + |\int_{K_{r} \times Y}f(\bfv, g)d\lambda(\bfv)d\mu_{Y}(g)|\\
 &\ll& (1 + r^{-3} + r^{-d-k -\kappa k +1}\|\frac{d\lambda}{d\bfx}\|_{\cC^{k}})\|f\|_{\cC^{1}(\bbS^{d-1}_{1}, \cS)}\zeta(\bdalp, e^{\frac{t}{2}})^{-\delta_{1}} + r^{-1}\|f\|_{\cC(\bbS_{1}^{d-1}, \cS)}e^{- \frac{d}{d-1}t}+ \|f\|_{L^{\infty}}r.
\eeqna
We then prove \eqref{est.2.10.1} by taking $r = \max\{\zeta(\bdalp, e^{\frac{t}{2}})^{-\delta_{2}}, e^{- \frac{d}{2(d-1)}t}\}$ with $\delta_{2} = \frac{1}{d+k +\kappa k }\delta_{1}$.
The proof of \eqref{est.2.10.2} follows that of \eqref{est.2.10.1}, with \eqref{est.3.2.2} in place of \eqref{est.3.2.1}, so that \eqref{est.4.3.1} is replaced by
\beqnas
& &|\int_{\bbS^{d-1}_{1}}f(\bfv, ME_{1}(\bfv)A^{t})d\lambda(\bfv) - \int_{\bbS^{d-1}_{1} \times X_{q}}f(\bfv, g)d\lambda(\bfv)d\mu_{q}(g)|\\
&\ll& (1 + r^{-3} + r^{-d-k-\kappa k+1}\|\frac{d\lambda}{d\bfx}\|_{\cC^{k}})\|f\|_{\cC^{1}(\bbS^{d-1}_{1}, \cS)}e^{- \lambda_{1} t} + r^{-1}\|f\|_{\cC(\bbS_{1}^{d-1}, \cS)}e^{- \frac{d}{d-1}t} + \|f\|_{L^{\infty}}r.
\eeqnas
Taking $r^{d+k+\kappa k} = e^{- \frac{\lambda_{1}}{2}t}$ (shrinking $\lambda_{1}$ if necessary so that $\lambda_{1} \le \frac{d(d+k+\kappa k)}{d-1}$), the right-hand side is $\ll \|f\|_{\cC^{1}(\bbS^{d-1}_{1}, \cS)}e^{- \lambda_{2} t}$ with $\lambda_{2} = \frac{\lambda_{1}}{2(d+k+\kappa k)}$, which proves \eqref{est.2.10.2}.
\epf

\section{Proof of Main Theorem}
\subsection{Estimates with respect to $\cE_{t}(\sigma)$}
Write 
\beqnas
\nonumber  \Theta_{t}(e^{\frac{t}{d-1}}\sigma)  
&=&   \{(\bfv, \bfx) \in \bbS^{d-1}_{1} \times \bbR^{d}, 0 \leq e^{-t}\|\bfx A^{-t}\| < \sigma^{d-1}, \bfx \in \bbR_{ \geq 0}\bfe_{1} + \cS^{d-1}_{1, \perp}\}.
\eeqnas
Noting that $e^{-t}\bfx A^{-t} = x_{1}\bfe_{1} + e^{-\frac{d}{d-1}t}\bfx_{\perp}$, we verify that the set $\Theta_{t}(e^{\frac{t}{d-1}}\sigma)$ (resp. $\Theta_{t, q}(e^{\frac{t}{d-1}}\sigma)$) is non-decreasing in $t$, such that $\Theta(\sigma) = \bigcup_{t} \Theta_{t}(e^{\frac{t}{d-1}}\sigma)$ .  Correspondingly, $\cE_{t}(\sigma)$ (resp. $\cE_{t, q}(\sigma)$) is non-increasing in $t$, yielding $\cE(\sigma) = \bigcap_{t} \cE_{t}(\sigma)$ and  $\cE_{q}(\sigma) = \bigcap_{t} \cE_{t, q}(\sigma)$. We remark that in contrast to the present setting where the monotonicity holds,  \cite{M-Sannals1} addresses a more general scenario by considering $\limsup$ and $\liminf$ of the respective sets.

We will make use of the following facts on the set $\cE_{t}(\sigma)$ in the proof of Theorem 1.3. 

\blem
There exists a function $S_{t}: Y \rightarrow \bbR^{+}$ such that
\beqna\label{st.eq}
\un_{\cE_{t}(\sigma)}(\bfv, g) =  \un_{\bbS^{d-1}_{1} }(\bfv)\un_{[\sigma^{d-1}, +\infty)}(S_{t}(g)).
\eeqna
Similarly, for each $q$ there exists a function $S_{t, q}: X_{q} \rightarrow \bbR^{+}$ satisfying 
\beqna
\label{st.eq.q}\un_{\cE_{t, q}(\sigma)}(\bfv, g) =  \un_{\bbS^{d-1}_{1} }(\bfv)\un_{[\sigma^{d-1}, +\infty)}(S_{t, q}(g)).
\eeqna

\elem

\bpf
By the definition~\eqref{def.ce.t}, we have
\beqnas
\cE_{t}(\sigma) = \{(\bfv, g) \in \bbS^{d-1}_{1} \times Y,  \hat{\cX}_{(t, \sigma)}(g) = 0\}, 
\eeqnas
where $\hat{\cX}_{(t, \sigma)}$ is a function on $Y$, defined by
$$
\hat{\cX}_{(t, \sigma)}(g) =  \sum_{\bfx \in \bbZ^{d}g}\un_{[0, \sigma^{d-1})}(e^{-t}\|\bfx A^{-t}\|)\un_{\bbR_{\geq 0}}(x_{1})\un_{\cS^{d-1}_{1, \perp}}(\bfx_{\perp}),
$$ 
and $\bfx_{\perp} = (0, x_{2}, \dots, x_{d})$ for $\bfx \in \bbR^{d}$. 
%
Define a function $s_{t}: \bbR^{d} \rightarrow \bbR^{+}$ as
\beqnas\label{def.st.1}
s_{t}(\bfx) = \left\{
\begin{array}{ll}
e^{-t}\|\bfx A^{-t}\|, & x_{1} \geq 0, \bfx_{\perp} \in \cS^{d-1}_{1, \perp}\\
+\infty, & otherwise,
\end{array}
\right.
\eeqnas
and a function $S_{t}: Y \rightarrow \bbR^{+}$ by
\beqna\label{def.st}
S_{t}(g) = \inf_{\bfx \in \bbZ^{d}g}s_{t}(\bfx),
\eeqna
such that  $\hat{\cX}_{(t, s)}(g)  = 0$ holds for any $s^{d-1} \in (0, S_{t}(g)]$.  Notice that the case $S_{t}(g) = +\infty$ corresponds to the situation where the ray in the direction $\bfe_{1}$ does not intersect any ball with radius $e^{-\frac{t}{d-1}}$ on the points of the affine lattice $g$. Thus, we derive that
 \beqnas
\un_{\cE_{t}(\sigma)}(\bfv, g) = \un_{\bbS^{d-1}_{1} }(\bfv)\un_{\{\hat{\cX}_{(t, \sigma)}(g) = 0\}}(g) =  \un_{\bbS^{d-1}_{1}}(\bfv)\un_{[\sigma^{d-1}, +\infty)}(S_{t}(g)).
\eeqnas
For $\cE_{t, q}(\sigma)$, define $S_{t, q}: X_{q} \rightarrow \bbR^{+}$ as
$$
S_{t, q}(g) = \inf_{\bfx \in (\bbZ^{d} + \frac{\bfp}{q})g}s_{t}(\bfx).
$$
The proof of \eqref{st.eq.q} is omitted, as it follows the same argument as $\cE_{t}(\sigma)$.
\epf
%
%
\blem
For $0 < \sigma_{1} \leq \sigma_{2}$ and $e^{-\frac{d}{d-1}t_{0}} < \sigma_{1}^{d-1}$, we have
\beqna
\label{est.4.2}| \int_{\bbS^{d-1}_{1}} (\un_{\cE_{t_{0}}(\sigma_{1})} -  \un_{\cE_{t_{0}}(\sigma_{2})}) (\bfv, (\un_{d}, \bdalp)ME_{1}(\bfv)A^{t}) d\lambda(\bfv)| \ll \sigma_{2}^{d-1} - \sigma_{1}^{d-1} + e^{-\frac{d}{d-1}t_{0}}, \\
\label{est.4.2.3}| \int_{\bbS^{d-1}_{1}} (\un_{\cE_{t_{0}, q}(\sigma_{1})} -  \un_{\cE_{t_{0}, q}(\sigma_{2})}) (\bfv, ME_{1}(\bfv)A^{t}) d\lambda(\bfv)| \ll \sigma_{2}^{d-1} - \sigma_{1}^{d-1} + e^{-\frac{d}{d-1}t_{0}}.
\eeqna
Moreover, 
\beqna
\label{est.4.2.2}| \int_{\bbS^{d-1}_{1}} (\un_{\cE(\sigma_{1})} -  \un_{\cE(\sigma_{2})}) (\bfv, (\un_{d}, \bdalp)ME_{1}(\bfv)A^{t}) d\lambda(\bfv)| \ll \sigma_{2}^{d-1} - \sigma_{1}^{d-1}, \\
\label{est.4.2.4}| \int_{\bbS^{d-1}_{1}} (\un_{\cE_{q}(\sigma_{1})} -  \un_{\cE_{q}(\sigma_{2})}) (\bfv, ME_{1}(\bfv)A^{t}) d\lambda(\bfv)| \ll \sigma_{2}^{d-1} - \sigma_{1}^{d-1}.
\eeqna
\elem
\bpf
By \eqref{st.eq}, we have
\beqnas
(\un_{\cE_{t_{0}}(\sigma_{1})} -  \un_{\cE_{t_{0}}(\sigma_{2})}) (\bfv, g) =  \un_{\bbS^{d-1}_{1} }(\bfv)\un_{[\sigma_{1}^{d-1}, \sigma^{d-1}_{2})}(S_{t_{0}}(g)),
\eeqnas
and $S_{t_{0}}\big((\un_{d}, \bdalp) ME_{1}(\bfv)A^{t} \big) \in [\sigma_{1}^{d-1}, \sigma^{d-1}_{2})$ implies that $e^{-t_{0}}\|(\bfz + \bdalp)ME_{1}(\bfv)A^{t}A^{-t_{0}}\|  \in  [\sigma_{1}^{d-1}, \sigma_{2}^{d-1})$ and $\|((\bfz + \bdalp)ME_{1}(\bfv)A^{t})_{\perp}\| < 1$ for some $\bfz \in \bbZ^{d}$. This leads to
\beqnas
e^{-t_{0}}\|((\bfz + \bdalp)ME_{1}(\bfv)A^{t}A^{-t_{0}})_{1}\| = |((\bfz + \bdalp)ME_{1}(\bfv)A^{t})_{1}| 
\geq \sigma_{1}^{d-1} - e^{-\frac{d}{d-1}t_{0}}.
\eeqnas
Since $E_{1}(\bfv) \in SO(d)$ keeps the volume of a region in $\bbR^{d}$, we derive that
\beqnas
\nonumber & &| \int_{\bbS^{d-1}_{1}} (\un_{\cE_{t_{0}}(\sigma_{1})} -  \un_{\cE_{t_{0}}(\sigma_{2})}) (\bfv, (\un_{d}, \bdalp)ME_{1}(\bfv)A^{t}) d\lambda(\bfv)|\\
\nonumber &\leq& \sum_{\bfz \in \bbZ^{d}}\int_{\bbS^{d-1}_{1}}\un_{[0, 1)}(e^{\frac{t}{d-1}} \|((\bfz + \bdalp)ME_{1}(\bfv))_{\perp}\|)\un_{[\sigma_{1}^{d-1} - e^{-\frac{d}{d-1}t_{0}}, \sigma_{2}^{d-1})}(e^{-t}|((\bfz + \bdalp)ME_{1}(\bfv))_{1}|)d\lambda(\bfv) \\
&\ll& (e^{-\frac{t}{d-1}})^{d-1}e^{t}(\sigma_{2}^{d-1} - (\sigma_{1}^{d-1}  - e^{-\frac{d}{d-1}t_{0}})) \ll \sigma_{2}^{d-1} - \sigma_{1}^{d-1} + e^{-\frac{d}{d-1}t_{0}}.
\eeqnas
where the implied constant depends on $(M, \bdalp)$ and $d$. Letting $t_{0} \rightarrow +\infty$ on both sides of \eqref{est.4.2}, we obtain \eqref{est.4.2.2}. 

The proof of \eqref{est.4.2.3} and \eqref{est.4.2.4} follows the same argument as $\cE_{t}(\sigma)$, $\cE(\sigma)$.
\epf

%
%
%

We recall the differentiability properties of the limiting distribution functions (Proposition 8.11, 8.13,  \cite{M-Sannals1}), which we reformulate as follows. 
\bprop\label{diffe}
Define $F_{q}(\sigma)$ for $\bdalp = \frac{\bfp}{q}$ and $F(\sigma)$  for irrational $\bdalp$ by
\beqnas
F_{q}(\sigma) &=& (\lambda \times \mu_{q})\{(\bfv, g) \in \bbS^{d-1}_{1} \times X_{q}, (\bbZ^{d}+  \frac{\bfp}{q})g \cap \Theta(\sigma)|_{\bfv} = \emptyset \},\\
F(\sigma) &=& (\lambda \times \mu_{Y})\{(\bfv, g) \in \bbS^{d-1}_{1} \times Y, \bbZ^{d}g \cap \Theta(\sigma)|_{\bfv} = \emptyset \}.
\eeqnas
Then $F_{q} \in C^{1}(0, +\infty)$ and  $F(\sigma) \in C^{2}(0, +\infty)$. Moreover, letting $\xi = \sigma^{d-1}$, $\frac{d}{d\xi}F_{q}(\sigma)$,  $\frac{d}{d\xi}F(\sigma)$ are uniformly bounded on $(0, +\infty)$. 
\eprop

\blem
For $0 < \sigma_{1} \leq \sigma_{2}$ and $e^{-\frac{d}{d-1}t_{0}} < \sigma_{1}^{d-1}$, we have
\beqna
\label{est.4.4.1}| \int_{\bbS^{d-1}_{1}}\int_{Y} (\un_{\cE(\sigma_{1})} -  \un_{\cE(\sigma_{2})}) (\bfv, g)  d\mu_{Y}(g)d\lambda(\bfv)| &\ll& |\sigma^{d-1}_{1} - \sigma^{d-1}_{2}|,\\
\label{est.4.4.2}| \int_{\bbS^{d-1}_{1}}\int_{Y} (\un_{\cE_{t_{0}}(\sigma)} -  \un_{\cE(\sigma)}) (\bfv, g)  d\mu_{Y}(g)d\lambda(\bfv)| &\ll& e^{-\frac{d}{d-1}t_{0}}.
\eeqna
Analogously, 
\beqna
\label{est.4.4.3}| \int_{\bbS^{d-1}_{1}}\int_{Y} (\un_{\cE_{q}(\sigma_{1})} -  \un_{\cE_{q}(\sigma_{2})}) (\bfv, g)  d\mu_{q}(g)d\lambda(\bfv)| &\ll& |\sigma^{d-1}_{1} - \sigma^{d-1}_{2}|,\\
\label{est.4.4.4}| \int_{\bbS^{d-1}_{1}}\int_{Y} (\un_{\cE_{t_{0}, q}(\sigma)} -  \un_{\cE_{q}(\sigma)}) (\bfv, g)  d\mu_{q}(g)d\lambda(\bfv)| &\ll& e^{-\frac{d}{d-1}t_{0}}.
\eeqna
\elem

\bpf
By Proposition~\ref{diffe}, we derive that
\beqnas
|\int_{\bbS^{d-1}_{1}}\int_{Y} (\un_{\cE(\sigma_{1})} -  \un_{\cE(\sigma_{2})}) (\bfv, g)  d\mu_{Y}(g)d\lambda(\bfv)| = |F(\sigma_{1}) - F(\sigma_{2})| \leq \max_{[\sigma^{d-1}_{1}, \sigma^{d-1}_{2}]}|F'(\xi)|\cdot |\sigma^{d-1}_{1} - \sigma^{d-1}_{2}|. 
\eeqnas
Observe that
\beqnas
\Theta(\sigma_{t_{0}}) = \{(\bfv, \bfx) \in \bbS^{d-1}_{1} \times \bbR^{d}, 0 \le x_{1} < \sigma^{d-1} - e^{-\frac{d}{d-1}t_{0}}, \bfx \in \bbR_{ \geq 0}\bfe_{1} + \cS^{d-1}_{1, \perp}\} \subset \Theta_{t_{0}}(e^{\frac{t_{0}}{d-1}}\sigma),
\eeqnas
where $\sigma^{d-1}_{t_{0}}=  \sigma^{d-1} - e^{-\frac{d}{d-1}t_{0}}$, which implies that $\cE_{t_{0}}(\sigma) \subset \cE(\sigma_{t_{0}})$. Since $\cE(\sigma) = \bigcap_{t} \cE_{t}(\sigma) \subset \cE_{t_{0}}(\sigma) $, we have
\beqnas
| \int_{\bbS^{d-1}_{1}}\int_{Y} (\un_{\cE_{t_{0}}(\sigma)} -  \un_{\cE(\sigma)}) (\bfv, g)  d\mu_{Y}(g)d\lambda(\bfv)| &\leq& | \int_{\bbS^{d-1}_{1}}\int_{Y} (\un_{\cE(\sigma_{t_{0}})} -  \un_{\cE(\sigma)}) (\bfv, g)  d\mu_{Y}(g)d\lambda(\bfv)|\\
&\ll& |\sigma^{d-1}_{t_{0}} - \sigma^{d-1}| = e^{-\frac{d}{d-1}t_{0}},
\eeqnas
where the last line is due to \eqref{est.4.4.1}. 
\epf

We introduce the following smooth approximations of $\un_{\cE_{t}(\sigma)}$. 

\blem\label{smt}
There exists a family of smooth, compactly supported functions $\{f_{\eps}\}$ on $ \bbS^{d-1}_{1} \times Y$ with $0 < \epsilon < \frac{1}{2}$, such that
\beqna\label{approxi}
0 \leq f_{\epsilon} \leq 1, \ \supp~f_{\epsilon} \subset \bbS^{d-1}_{1} \times \{y \in Y, \htt(\pi(y)) \leq 2\eps^{-\frac{1}{d}}\},
\eeqna
and
\beqna
\label{approxi.1}\|f_{\epsilon} - \un_{\cE_{t}(\sigma)}\|_{L^{1}(\bbS^{d-1}_{1} \times Y)} &\ll& \eps + e^{-\frac{2d}{d-1}t}, \\ 
\label{approxi.2}\|f_{\epsilon}\|_{\cC^{1}(\bbS^{d-1}_{1}, \cS)} &\ll& \epsilon^{- \frac{k(d+1)}{d}},  
\eeqna
\elem

\bpf
We construct the family $\{f_{\eps}\}_{\eps > 0}$ as follows. Since $\cE_{t}(\sigma)$ is unbounded in $\bbS^{d-1}_{1} \times Y$, we truncate it by defining
$$
\cE^{\eps}_{t}(\sigma) = \cE_{t}(\sigma) \cap (\bbS^{d-1}_{1} \times \{y \in Y, \htt(\pi(y)) \leq \eps^{-\frac{1}{d}}\}),
$$
so that
\beqna\label{eee}
\un_{\cE^{\eps}_{t}(\sigma)}(\bfv, g) = \un_{\cE_{t}(\sigma)}(\bfv, g)\un_{K(\eps^{-\frac{1}{d}})}(\pi(g)).
\eeqna
Then \eqref{est.ht.2} implies that 
$$
\mu_{Y}(\cE_{t}(\sigma)|_{\bfv} \cap \{\htt(\pi(y)) \geq \eps^{-\frac{1}{d}}\}) \ll \eps.
$$
Consequently, 
\beqna\label{est.4.5.11}
\|\un_{\cE^{\eps}_{t}(\sigma)} - \un_{\cE_{t}(\sigma)}\|_{L^{1}(\bbS^{d-1}_{1}\times Y)}\ll \eps. 
\eeqna
Let $K_{\eps} \in \cC^{\infty}_{c}(Y \times Y)$ be a mollifier kernel on $Y$, supported on $\{d^{Y}(g, g') < \eps\}$ and satisfying
\beqnas
K_{\eps} \geq 0, \qquad \|K_{\eps}(\cdot, g')\|_{\cC^{k}} \ll \eps^{-dim Y - k}, \qquad \int_{Y}K_{\eps}(g, g')d\mu_{Y}(g') = 1. 
\eeqnas
Define 
$$
f_{\eps}(\bfv, g) = \int_{Y} \un_{\cE^{\eps}_{t}(\sigma)}(\bfv, g') K_{\eps}(g, g') d\mu_{Y}(g').
$$
Then $0 \le f_{\eps} \le 1$, and $f_{\eps}(\bfv, \cdot)$ is supported on $\bbS^{d-1}_{1} \times \{y \in Y, \htt \cdot \pi(y) \leq 2\eps^{-\frac{1}{d}}\}$, since we can take $\eps$ sufficiently small so that $\htt \cdot \pi(g) \ll e^{d^{Y}(g, g')}\htt \cdot \pi(g') \leq 2\eps^{-\frac{1}{d}}$. For any differential operator $\cD$ of degree $k$, 
$$
\|\cD f_{\eps}(\bfv, \cdot)\|_{L^{\infty}} \leq \sup_{g \in K(2\eps^{-\frac{1}{d}})}\int_{B^{Y}(g, \eps)} |\cD K_{\eps}(g, g')| d\mu_{Y}(g') \ll \eps^{-k},
$$
using $\mu_{Y}(B^{Y}(g, \eps)) \asymp \eps^{dim Y}$. This yields
\beqnas
\|f_{\eps}\|_{\cC^{1}(S^{d-1}_{1}, \cS)} \ll \eps^{-\frac{k}{d}}\eps^{-k} = \eps^{-\frac{k(d+1)}{d}},
\eeqnas
which is \eqref{approxi.2}. Moreover, for each $\bfv$, 
\beqnas
\|f_{\eps}(\bfv, \cdot) - \un_{\cE^{\eps}_{t}(\sigma)}(\bfv, \cdot)\|_{L^{1}(Y)} \leq \int_{Y}\int_{Y}| \un_{\cE^{\eps}_{t}(\sigma)}(\bfv, g) -  \un_{\cE^{\eps}_{t}(\sigma)}(\bfv, g')| K_{\eps}(g, g')d\mu_{Y}(g')d\mu_{Y}(g). 
\eeqnas
Observe that the integral is non-trivial only when $d^{Y}(g,g') < \eps$. We claim that, whenever
$d^{Y}(g,g') < \eps$,
\beqna\label{flip}
|\un_{\cE^{\eps}_{t}(\sigma)}(\bfv, g) - \un_{\cE^{\eps}_{t}(\sigma)}(\bfv, g')|
\le \un_{\bbS^{d-1}_{1}}(\bfv)\big(\un_{W_{\eps}}(g) + \un_{T_{\eps}}(g)\big),
\eeqna
where $W_{\eps}$, $T_{\eps}$ are defined below. To prove \eqref{flip}, we write $T = \eps^{-1/d}$ and consider two cases.\\
\emph{Case 1:} $\un_{\cE_t}(\bfv,g) \neq \un_{\cE_t}(\bfv,g')$. Then there exists $\bfx' \in \bbZ^{d}g' \cap \Theta_t(e^{\frac{t}{d-1}}\sigma)$ while $\bbZ^{d}g \cap \Theta_t(e^{\frac{t}{d-1}}\sigma) = \emptyset$ (the other case is symmetric).
Since $\Theta_t(e^{\frac{t}{d-1}}\sigma) \subset B^{\bbR^{d}}(0, \sigma^{d-1}+1)$, we have $\|\bfx'\| \le \sigma^{d-1}+1$.  Writing $g = g'h$ with $h = (A, \bfz) \in Y$ and $d^{Y}(id, h) = d^{Y}(g, g') < \eps$, we obtain that $\bfx = \bfx' h$ satisfies
\beqnas
\|\bfx - \bfx'\| &=& \|\bfx'(I-A) - z\| \le (1+\|\bfx'\|)(\|I-A\|+\|z\|) \\
&\leq& (1+ \|\bfx'\|)\|I - h\|_{op} \ll  (\sigma^{d-1} +2) d^{Y}(\id, h) \ll (\sigma^{d-1} + 2)\eps. 
\eeqnas
Thus, the point $\bfx$ lies in the $C(d,\sigma)\eps$-thickening of $\partial\Theta_t(e^{\frac{t}{d-1}}\sigma)$ inside $B(0, \sigma^{d-1}+ 1)$, where $\partial\Theta_t(e^{\frac{t}{d-1}}\sigma) = \{x_1=0\} \bigcup \{x_1^{2} + e^{-\frac{2d}{d-1}t}\|\bfx_{\perp}\|^{2} = \sigma^{2(d-1)}, \|\bfx_{\perp}\| \le 1\} \bigcup \{\|\bfx_{\perp}\|=1, 0 < x_{1} < \sqrt{\sigma^{2(d-1)} - e^{-\frac{2d}{d-1}t}}\}$. This implies $g \in W_{\eps}$, where $W_{\eps}$ is given by
\beqnas
W_{\eps} = \{g \in Y, \bbZ^{d}g \cap A_{\eps} \neq \emptyset \}, 
\eeqnas
and
\beqnas
 A_{\eps} &=& B^{\bbR^{d}}(0, \sigma^{d-1}+1) \cap \big( \{\bfx \in \bbR^{d}, \sigma^{d-1} - C(\sigma, d)(\eps + e^{-\frac{2d}{d-1}t}) \leq x_{1} \leq \sigma^{d-1} + C(\sigma, d)\eps, \|\bfx_{\perp}\| \leq 1\} \\
 & &\cup \{0 < x_{1} <  \sigma^{d-1}+ C(\sigma, d)\eps,  \ 1 - C(\sigma, d)\eps \leq \|\bfx_{\perp}\| \leq 1 + C(\sigma, d)\eps\} \\
 & &\cup \{|x_{1}| \leq C(\sigma, d)\eps,  \|\bfx_{\perp}\| \leq 1 + C(\sigma, d)\eps\} \big),
\eeqnas
with $\vol(A_{\eps}) \ll  C(\sigma, d)(\eps + e^{-\frac{2d}{d-1}t}) $.\\
%
\emph{Case 2:} $\un_{\cE_t}(\bfv,g) = \un_{\cE_t}(\bfv,g')$. Then $\{\un_{\cE^{\eps}_{t}(\sigma)}(\bfv, g) \neq \un_{\cE^{\eps}_{t}(\sigma)}(\bfv, g')\}$ only comes from the truncation: either $\htt(\pi (g)) \le T < \htt(\pi(g'))$ or $\htt(\pi(g')) \le T < \htt(\pi (g))$. Let $\st(g) = (\htt(\pi (g)))^{-1}$.  
If $\htt(\pi(g)) \le T < \htt(\pi(g'))$, then $\st(g') < \eps^{1/d}$ and
$\st(g) \le e^{C\eps}(\st(g')+C\eps) \le e^{C\eps}(\eps^{1/d}+C\eps)$, hence
$T \ge \htt(\pi(g)) \ge T(1 - C'\eps^{\frac{d-1}{d}})$.
If $\htt(\pi(g')) \le T < \htt(\pi(g))$, then $\st(g) \ge e^{-C\eps}\st(g') - C\eps \ge e^{-C\eps}\eps^{1/d} - C\eps
\ge \eps^{1/d}(1 - C'\eps^{\frac{d-1}{d}})$, hence $T \le \htt(\pi(g)) \le T(1 + C'\eps^{\frac{d-1}{d}})$.
In both cases
\beqnas 
g \in T_{\eps} := \{\htt \cdot \pi \in  (T(1 - C'\eps^{\frac{d-1}{d}}),  T(1 + C'\eps^{\frac{d-1}{d}})) \}.
\eeqnas
This proves \eqref{flip}. 

By Siegel's formula, 
\beqnas
\int_{Y}\sum_{\bfx \in \bbZ^{d}g}\un_{A_{\eps}}(\bfx)d\mu_{Y}(g) =  \vol(A_{\eps}),
\eeqnas
so that
\beqnas
\mu_{Y}(W_{\eps}) \leq \int_{Y}\sum_{\bfx \in \bbZ^{d}g}\un_{A_{\eps}}(\bfx)d\mu_{Y}(g) = \vol(A_{\eps}) \leq  C(\sigma, d)(\eps + e^{-\frac{2d}{d-1}t}) .
\eeqnas
Also, by  \eqref{est.ht.2}, we have $\mu_{Y}(T_{\eps}) \ll \eps$. Therefore, 
\beqnas
\|f_{\eps} - \un_{\cE^{\eps}_{t}(\sigma)}\|_{L^{1}(\bbS^{d-1}_{1} \times Y)} &\leq& \int_{Y}(\un_{W_{\eps}}(g) + \un_{T_{\eps}}(g))(\int_{Y}K_{\eps}(g, g')d\mu_{Y}(g'))d\mu_{Y}(g) \\
&=& \mu_{Y}(W_{\eps}) +  \mu_{Y}(T_{\eps}) \ll  C(\sigma, d)(\eps + e^{-\frac{2d}{d-1}t}). 
\eeqnas
Together with \eqref{est.4.5.11}, we derive \eqref{approxi.1}.
\epf

%
%

\subsection{Proof of Theorem 1.3}

It suffices to prove the irrational case; the rational case follows similarly. We first extend the estimates in Proposition~\ref{prop2.10} to the indicator function of $\cE_{t_{0}}(\sigma)$ and $\cE(\sigma)$.
\blem\label{l3.4}
Let $\lambda$ be a probability measure on $\bbS^{d-1}_{1}$ with smooth density with respect to the Lebesgue measure. For $\boldsymbol{\alpha}  \in \bbR^{d} \setminus \bbQ^{d}$, there exist positive constants $\bar{\delta}$, $\tilde{\kappa}$, such that
\beqna\label{est.4.5.1}
\nonumber |\int_{\bbS^{d-1}_{1}} \un_{\cE_{t_{0}}(\sigma)} (\bfv, (\un_{d}, \bdalp)ME_{1}(\bfv)A^{t}) d\lambda(\bfv) - \int_{\bbS^{d-1}_{1} \times Y}\un_{\cE_{t_{0}}(\sigma)}(\bfv, g) d\lambda(\bfv)d\mu_{Y}(g)|  \\
 \ll\zeta(\bdalp, e^{\frac{t}{2}})^{-\bar{\delta}} + e^{- \tilde{\kappa} t} + e^{-\frac{2d}{d-1}t_{0}}.
\eeqna
Moreover, 
\beqna\label{est.4.5.3}
\nonumber |\int_{\bbS^{d-1}_{1}} \un_{\cE(\sigma)} (\bfv, (\un_{d}, \bdalp) ME_{1}(\bfv)A^{t}) d\lambda(\bfv) - \int_{\bbS^{d-1}_{1} \times Y}\un_{\cE(\sigma)}(\bfv, g) d\lambda(\bfv)d\mu_{Y}(g)| \\
 \ll  \zeta(\bdalp, e^{\frac{t}{2}})^{-\bar{\delta}} + e^{- \tilde{\kappa} t}.
\eeqna
 The implied constants depend on $M$, $\boldsymbol{\alpha}$,  $d$, $\lambda$. 
\elem

\bpf


Let $\{f_{\epsilon}\}$ be the family of smooth approximations to $\un_{\cE_{t_{0}}(\sigma)}$ in Lemma~\ref{smt}. We separate the difference into three parts,
\beqnas
& &|\int_{\bbS^{d-1}_{1}}  \un_{\cE_{t_{0}}(\sigma)}(\bfv, (\un_{d}, \bdalp) ME_{1}(\bfv)A^{t}) d\lambda(\bfv) - \int_{\bbS^{d-1}_{1} \times Y}  \un_{\cE_{t_{0}}(\sigma)}(\bfv, g) d\mu_{Y}(g)d\lambda(\bfv)|\\
&\leq& |\int_{\bbS^{d-1}_{1}}   f_{\epsilon}(\bfv, (\un_{d}, \bdalp) ME_{1}(\bfv)A^{t}) d\lambda(\bfv) - \int_{\bbS^{d-1}_{1} \times Y} f_{\epsilon}(\bfv, g) d\mu_{Y}(g)d\lambda(\bfv)| \\
& & + |\int_{\bbS^{d-1}_{1}}  \un_{\cE_{t_{0}}(\sigma)}(\bfv, (\un_{d}, \bdalp) ME_{1}(\bfv)A^{t}) d\lambda(\bfv) - \int_{\bbS^{d-1}_{1}} f_{\epsilon}(\bfv, (\un_{d}, \bdalp) ME_{1}(\bfv)A^{t}) d\lambda(\bfv)| \\
& &+ |\int_{\bbS^{d-1}_{1} \times Y}  \un_{\cE_{t_{0}}(\sigma)}(\bfv, g) d\mu_{Y}(g)d\lambda(\bfv) -\int_{\bbS^{d-1}_{1} \times Y} f_{\epsilon}(\bfv, g) d\mu_{Y}(g)d\lambda(\bfv)|\\
&:=& I_{1} + I_{2} + I_{3}.
\eeqnas
For the first term $I_{1}$, applying \eqref{est.2.10.1} to $f_{\epsilon}$ we derive
\beqna\label{est.1}
\nonumber I_{1} &\ll&  \|f_{\epsilon}\|_{\cC^{1}(\bbS^{d-1}_{1}, \cS)}(\zeta(\bdalp, e^{\frac{t}{2}})^{-\delta_{2}} + e^{- \frac{d}{2(d-1)} t}) \\
&\ll&   \epsilon^{- \frac{k(d+1)}{d}}(\zeta(\bdalp, e^{\frac{t}{2}})^{-\delta_{2}} + e^{- \frac{d}{2(d-1)} t}),
\eeqna
where the last line is due to \eqref{approxi.2}. 
Note that, since $\cE^{\eps}_{t_{0}}(\sigma) = \cE_{t_{0}}(\sigma) \cap \pi^{-1}(K(\eps^{-\frac{1}{d}}))$ and $\int_{Y}K_{\eps}(\cdot, g')d\mu_{Y}(g') = 1$, we have by \eqref{flip}, for $ (\bfv, g) \in \bbS^{d-1}_{1} \times Y$, 
$$
|\un_{\cE_{t_{0}}(\sigma)}(\bfv, g) - f_{\epsilon}(\bfv, g)| \le \un_{K^{c}(\eps^{-\frac{1}{d}})}(\pi(g)) + \un_{W_{\eps}}(g) + \un_{T_{\eps}}(g).
$$
Hence, writing $g_{t}(\bfv) = (\un_{d}, \bdalp)ME_{1}(\bfv)A^{t}$,
\beqna\label{est.2}
\nonumber I_{2} &\leq& \int_{\bbS^{d-1}_{1}}\un_{K^{c}(\eps^{-\frac{1}{d}})}(ME_{1}(\bfv)A^{t}) )d\lambda(\bfv) + \int_{\bbS^{d-1}_{1}}\un_{W_{\eps}}(g_{t}(\bfv))d\lambda(\bfv) + \int_{\bbS^{d-1}_{1}}\un_{T_{\eps}}(g_{t}(\bfv))d\lambda(\bfv) \ll \eps + e^{-\frac{2d}{d-1}t_{0}} ,
\eeqna
where each term is bounded by the same lattice-point estimate as in the proof of \eqref{est.4.2}: 
\beqnas
 |\int_{\bbS^{d-1}_{1}}\un_{K^{c}(\eps^{-\frac{1}{d}})}(ME_{1}(\bfv)A^{t})d\lambda(\bfv)| &\le&  \lambda\{\bfv: \htt(\pi(g_{t}(\bfv))) \ge \eps^{-\frac{1}{d}}\} \asymp \eps, \\
\int_{\bbS^{d-1}_{1}}\un_{W_{\eps}}(g_{t}(\bfv))d\lambda(\bfv) &\le& \sum_{\bfz \in \bbZ^{d}}\int_{\bbS^{d-1}_{1}}\un_{A_{\eps}}((\bfz + \bdalp)ME_{1}(\bfv)A^{t})d\lambda(\bfv) \ll \vol(A_{\eps}) \ll \eps +  e^{-\frac{2d}{d-1}t_{0}}, \\
\int_{\bbS^{d-1}_{1}}\un_{T_{\eps}}(g_{t}(\bfv))d\lambda(\bfv) &\le& \lambda\{\bfv: \htt(\pi(g_{t}(\bfv))) \ge T(1-C'\eps^{\frac{d-1}{d}})\} \ll (T(1-C'\eps^{\frac{d-1}{d}}))^{-d} \asymp \eps,
\eeqnas 
since $E_{1}(\bfv)$, $A^{t}$ and $M$ are all unimodular. By \eqref{approxi.1} we have 
\beqna\label{est.3}
I_{3} \leq \epsilon  + e^{-\frac{2d}{d-1}t_{0}} .
\eeqna
Combining \eqref{est.1}, \eqref{est.3} and \eqref{est.2}, we prove \eqref{est.4.5.1} by taking $\epsilon = (\zeta(\bdalp, e^{\frac{t}{2}})^{-\delta_{2}} + e^{- \frac{d}{2(d-1)} t})^{\frac{d}{k(d+1)+d}}$, $\bar{\delta} = \frac{d\delta_{2}}{k(d+1)+d}$ and $\tilde{\kappa} = \frac{d^{2}}{2(d-1)(k(d+1)+d)}$. We derive \eqref{est.4.5.3} by  letting $t_{0} \rightarrow +\infty$ on the right side of \eqref{est.4.5.1}. 

\epf

We now proceed to prove Theorem~\ref{m.1}.

\bpf
Notice that for $t > t_{0}$, we have $\cE(\sigma) \subset \cE_{t}(\sigma) \subset \cE_{t_{0}}(\sigma)$, so that
\beqna\label{bdd.lu}
\nonumber \int_{\bbS^{d-1}_{1}} \un_{\cE(\sigma)} (\bfv, (\un_{d}, \bdalp)ME_{1}(\bfv)A^{t}) d\lambda(\bfv) &\leq& \int_{\bbS^{d-1}_{1}} \un_{\cE_{t}(\sigma)} (\bfv,  (\un_{d}, \bdalp)ME_{1}(\bfv)A^{t}) d\lambda(\bfv)\\
& \leq& \int_{\bbS^{d-1}_{1}} \un_{\cE_{t_{0}}(\sigma)} (\bfv, (\un_{d}, \bdalp)ME_{1}(\bfv)A^{t}) d\lambda(\bfv).
\eeqna
Applying \eqref{est.4.5.1}, we obtain
\beqna\label{ineq.u}
\nonumber & & \int_{\bbS^{d-1}_{1}}  \un_{\cE_{t_{0}}(\sigma)} (\bfv, (\un_{d}, \bdalp)ME_{1}(\bfv)A^{t}) d\lambda(\bfv)  - \int_{\bbS^{d-1}_{1}} \int_{Y} \un_{\cE(\sigma)}(\bfv, g) d\mu_{Y}(g)d\lambda(\bfv) \\
\nonumber &\ll& \int_{\bbS^{d-1}_{1} \times Y} \un_{\cE_{t_{0}}(\sigma)} (\bfv, g) d\mu_{Y}(g)d\lambda(\bfv) - \int_{\bbS^{d-1}_{1}} \int_{Y} \un_{\cE(\sigma)}(\bfv, g) d\mu_{Y}(g)d\lambda(\bfv)\\
 \nonumber& &+ \zeta(\bdalp, e^{\frac{t}{2}})^{-\bar{\delta}} + e^{-  \tilde{\kappa} t}+ e^{- \frac{2d}{d-1} t_{0}}\\
 &\ll& \zeta(\bdalp, e^{\frac{t}{2}})^{-\bar{\delta}} + e^{- \tilde{\kappa} t} + e^{- \frac{2d}{d-1} t_{0}},
\eeqna
where the last line is due to \eqref{est.4.4.2}. By \eqref{est.4.5.3}, 
\beqna\label{ineq.l}
\nonumber \int_{\bbS^{d-1}_{1}} \un_{\cE(\sigma)} (\bfv, (\un_{d}, \bdalp)ME_{1}(\bfv)A^{t}) d\lambda(\bfv) - \int_{\bbS^{d-1}_{1}} \int_{Y} \un_{\cE(\sigma)}(\bfv, g) d\mu_{Y}(g)d\lambda(\bfv) \\
\gg -  \big(\zeta(\bdalp, e^{\frac{t}{2}})^{-\bar{\delta}} + e^{-\tilde{\kappa} t} \big).
\eeqna
Combining \eqref{bdd.lu}, \eqref{ineq.u} and \eqref{ineq.l} and taking $t_{0} \rightarrow t^{-}$, we have
\beqna\label{ineq.1.5}
\nonumber |\int_{\bbS^{d-1}_{1}} \un_{\cE_{t}(\sigma)}(\bfv, (\un_{d}, \bdalp)ME_{1}(\bfv)A^{t}) d\lambda(\bfv)  - \int_{\bbS^{d-1}_{1} \times Y} \un_{\cE(\sigma)}(\bfv, g)d\lambda(\bfv) d\mu_{Y}(g)|\\
 \ll   \zeta(\bdalp, e^{\frac{t}{2}})^{-\bar{\delta}} + e^{-  \tilde{\kappa} t} + e^{- \frac{2d}{d-1} t}.
 \eeqna
By \eqref{key.1}, we have
\beqnas\label{ineq.1.6}
& &| \lambda(\bfv, r^{d-1}\tau_{r}(\bfq, \bfv) \geq \xi) - \int^{\infty}_{\xi}\Phi(x)dx| \\
&\leq& |\int_{\bbS^{d-1}_{1}} \un_{\cE_{t}(\sigma)}(\bfv, (\un_{d}, \bdalp)ME_{1}(\bfv)A^{t}) d\lambda(\bfv)  - \int_{\bbS^{d-1}_{1} \times Y} \un_{\cE(\sigma)}(\bfv, g)d\lambda(\bfv) d\mu_{Y}(g)|\\
& & \qquad +  \max\{\epsilon_{1, t}, \epsilon_{2, t}\},
\eeqnas
where
\beqnas
\epsilon_{1, t} &=& | \int_{\bbS^{d-1}_{1}}(\un_{\cE_{t}(\bar{\sigma}_{t})} - \un_{\cE_{t}(\sigma)})(\bfv, (\un_{d}, \bdalp)ME_{1}(\bfv)A^{t}) d\lambda(\bfv)|, \\
 \epsilon_{2, t} &=& | \int_{\bbS^{d-1}_{1}}(\un_{\cE_{t}(\tilde{\sigma}_{t})} - \un_{\cE_{t}(\sigma)})(\bfv,  (\un_{d}, \bdalp)ME_{1}(\bfv)A^{t}) d\lambda(\bfv)|.
\eeqnas
By \eqref{est.4.2}, we derive that
\beqna\label{est.err}
\epsilon_{1, t} \ll |\bar{\sigma}^{d-1}_{t} - \sigma^{d-1}| + e^{-\frac{d}{d-1}t} \ll e^{-\frac{d}{d-1}t}, \ \epsilon_{2, t}  \ll |\tilde{\sigma}^{d-1}_{t} - \sigma^{d-1}| + e^{-\frac{d}{d-1}t} \ll e^{-\frac{d}{d-1}t}.
\eeqna
Thus combining \eqref{ineq.1.5} and \eqref{est.err}, we obtain
\beqnas
| \lambda(\bfv, r^{d-1}\tau_{r}(\bfq, \bfv) \geq \xi) - \int^{\infty}_{\xi}\Phi(x)dx| \ll \zeta(\bdalp, e^{\frac{t}{2}})^{-\bar{\delta}} + e^{-  \min\{ \tilde{\kappa}, \frac{d}{d-1}\} t}.
\eeqnas
We complete the proof of \eqref{est.main.l} by noticing $r^{d-1}= e^{-t}$ and letting $\bar{\kappa} =  \min\{(d-1) \tilde{\kappa}, d\}$.

\epf

\subsection{Proof of Corollary~\ref{c.1}}
We first recall some facts on Diophantine type vectors from Section 2.1.2, \cite{L26}. A vector $\bfb \in \bbR^{d}$ is said to be of Diophantine type $\kappa \geq 1$ if there exists $c_{\kappa} > 0$ such that
\beqna\label{da.k}
|\bfb - \frac{\bfp}{q}| > c_{\kappa}q^{-\kappa}
\eeqna
for any $\bfp \in \bbZ^{d}$ and $q \in \bbN$, where $|\cdot|$ is the supremum norm of $\bbR^{d}$. Note that for a vector $\bfb$ of Diophantine type $\kappa$, one has
 \beqna\label{est.z}
\zeta(\bfb, T) \geq c^{\frac{1}{\kappa+1}}_{\kappa} T^{\frac{1}{\kappa+1}}.
\eeqna
For $c>0$, define the set
\beqnas
D(\kappa, c) = \{\bfb \in \bbT^{d} , |q\bfb|_{\bbZ} >  cq^{-\kappa+1}, \ \forall q \in \bbN\},
\eeqnas
which is decreasing in $c$. 
For $\kappa > d+1$ and $\epsilon > 0$, the set $\bigcup^{\infty}_{c \geq \epsilon}D(\kappa, c)$ consists of vectors of Diophantine type $\kappa$, and 
\beqna\label{est.m}
m_{\bbT^{d}}\big(\bbT^{d} \setminus \bigcup^{\infty}_{c \geq \epsilon} D(\kappa, c) \big) \ll \epsilon.
\eeqna



\noindent
\emph{Proof of Corollary~\ref{c.1}.}
Let $\cF$ be the fundamental parallelogram for $\cL$ and let $p$ be the probability density of $\Lambda$ with respect to the Lebesgue measure. Notice that
\beqnas
& &\Lambda((\bfQ, \bfV) \in T^{1}(r^{d-1}\cK_{r}): \cT_{r}(\bfQ, \bfV) \geq \xi ) \\
&=& \int_{\bbS^{d-1}_{1}}\int_{\bbR^{d}}\un_{\{r^{d-1}\tau_{r}(\bfq, \bfv) \geq \xi \}}p(r^{d-1}\bfq, \bfv)r^{d(d-1)}d\bfq d\bfv\\
&=& \int_{\bbS^{d-1}_{1}}\int_{\cF}\un_{\{r^{d-1}\tau_{r}(\bfq_{0}, v) \geq \xi \}}(\sum_{\bfq \in \bfq_{0} + \cL}p(r^{d-1}\bfq, \bfv)r^{d(d-1)})d\bfq_{0}d\bfv.
\eeqnas
For $\bfq_{0} \in \cF$ and $\bfv \in \bbS^{d-1}_{1}$, we define
\beqnas
F_{r}(\bfq_{0}, \bfv) := \sum_{\bfq \in \bfq_{0} + \cL}p(r^{d-1}\bfq, \bfv)r^{d(d-1)},  \quad  F(\bfv) = \int_{\bbR^{d}}p(\bfq, \bfv)d\bfq. 
\eeqnas
Then $F(\bfv)$ can be viewed as a probability density on $\bbS^{d-1}_{1}$, since $\int_{\bbS^{d-1}_{1}}F(\bfv)d\bfv = \int_{\bbS^{d-1}_{1}}\int_{\bbR^{d}}p(\bfq, \bfv)d\bfq d\bfv =1$. 
\blem
Under the assumption~\eqref{asp.c1}, it holds that
\beqna\label{est.fr}
\sup_{\bfq_{0} \in \cF, \bfv \in \bbS^{d-1}_{1}}|F_{r}(\bfq_{0}, \bfv) - F(\bfv) | \ll r^{(d-1)\gamma}.
\eeqna
\elem
\bpf
Let $\epsilon = r^{d-1}$ and $g(\bfq) := p(\epsilon \bfq, \bfv)\epsilon^{d}$, then $F_{r}(\bfq_{0}, \bfv) = \sum_{q \in q_{0} + \cL}g(\bfq)$. By the Poisson summation formula on the unimodular lattice, we obtain
\beqna\label{eq.ps}
\sum_{\bfq \in \bfq_{0} + \cL}g(\bfq) = \sum_{\bfk \in \cL^{*}}e^{2\pi i \bfq_{0} \cdot \bfk}\hat{g}(\bfk), 
\eeqna
where $\cL^{*} = \bbZ^{d}M^{-tr}$ is the dual lattice of $\cL$ and $\hat{g}$ is the Fourier transform of $g$, given by
\beqnas
\hat{g}(\bfk) = \int_{\bbR^{d}}p(\epsilon \bfx, \bfv)\epsilon^{d}e^{-2\pi i \bfx \cdot \bfk}d\bfx, 
= \hat{p}_{\bfv}(\frac{\bfk}{\epsilon}), \quad \hat{p}_{\bfv}(\bfy) = \int_{\bbR^{d}}p(\bfx, \bfv)e^{-2\pi i \bfx \cdot \bfy}d\bfx.
\eeqnas
Inserting it into \eqref{eq.ps}, we derive $F_{r}(\bfq_{0}, \bfv) = \sum_{\bfk \in \cL^{*}}e^{2\pi i \bfq_{0} \cdot \bfk}\hat{p}_{\bfv}(\frac{\bfk}{\epsilon})$. Notice that $\hat{p}_{\bfv}(\mathbf{0}) = F(\bfv)$, so that
\beqna\label{est.Fr}
F_{r}(\bfq_{0}, \bfv) - F(\bfv) = \sum_{\bfk \in \cL^{*}, \bfk \neq \mathbf{0}} e^{2\pi i \bfq_{0} \cdot \bfk}\hat{p}_{\bfv}(\frac{\bfk}{\epsilon}).
\eeqna
Recall that $\widehat{\partial^{\gamma}_{i}p_{\bfv}}(\bfm) = \int_{\bbR^{d}}\partial^{\gamma}_{q_{i}}p(\bfq, \bfv)e^{-2\pi i\bfq \cdot \bfm}d\bfq$ and $\widehat{\partial_{i}^{\gamma}p_{\bfv}}(\bfm) = (2\pi i m_{i})^{\gamma}\hat{p}_{\bfv}(\bfm)$ for $\gamma > d$, so that by the assumption \eqref{asp.c1}, for any $i=1, \dots, d$, 
\beqna\label{ineq.4.6}
|\hat{p}_{\bfv}(\bfm)| \ll |m_{i}^{-\gamma}| \int_{\bbR^{d}}|\partial^{\gamma}_{q_{i}}p(\bfq, \bfv)|d\bfq. 
\eeqna
Assume that $|m_{j}| = \max_{1 \leq i \leq d}|m_{i}|$, such that $|m_{j}| \geq \frac{1}{\sqrt{d}}\|\bfm\|$ and $|m_{j}^{-\gamma}| \ll \|\bfm\|^{-\gamma}$. Then taking $i=j$ in \eqref{ineq.4.6}, we derive
\beqnas
|\hat{p}_{\bfv}(\frac{\bfk}{\epsilon})|  \ll \epsilon^{\gamma}\|\bfk\|^{-\gamma}\sup_{\bfv \in \bbS^{d-1}_{1}}\sum^{d}_{i=1}\|\partial^{\gamma}_{i}p(\cdot, \bfv)\|_{L^{1}(\bbR^{d})},
\eeqnas
which leads to
\beqnas
 \sup_{\bfq_{0} \in \cF, \bfv \in \bbS^{d-1}_{1}}|F_{r}(\bfq_{0}, \bfv) - F(\bfv) | &\ll& \epsilon^{\gamma} \sup_{\bfv \in \bbS^{d-1}_{1}}\sum^{d}_{i=1}\|\partial^{\gamma}_{i}p(\cdot, \bfv)\|_{L^{1}(\bbR^{d})}\sum_{\bfk \in \cL^{*}, \bfk \neq \mathbf{0}}\|\bfk\|^{-\gamma} \\
 &\ll& r^{(d-1)\gamma}. 
\eeqnas
 by \eqref{est.Fr} and the assumption \eqref{asp.c1}. 
\epf
Then,  
\beqnas
& &|\Lambda((\bfQ, \bfV) \in T^{1}(r^{d-1}\cK_{r}): \cT_{r}(\bfQ, \bfV) \geq \xi )- \int^{+\infty}_{\xi}\Phi(x)dx| \\
&=& |\int_{\bbS^{d-1}_{1}}\int_{\cF}\un_{\{r^{d-1}\tau_{r}(\bfq_{0}, v) \geq \xi \}}F_{r}(\bfq_{0}, \bfv)d\bfq_{0}d\bfv - \int^{+\infty}_{\xi}\Phi(x)dx |\\
&\leq& \int_{\cF}|\int_{\bbS^{d-1}_{1}}\un_{\{r^{d-1}\tau_{r}(\bfq_{0}, v) \geq \xi\}}F(\bfv)d\bfv - \int^{+\infty}_{\xi}\Phi(x)dx| d\bfq_{0}  + \sup_{\bfq_{0} \in \cF, \bfv \in \bbS^{d-1}_{1}} |F_{r}(\bfq_{0}, \bfv) - F(\bfv)|.
\eeqnas
Notice that $\bdalp = -\bfq_{0}M^{-1}$, 
\beqnas
& & \int_{\cF}|\int_{\bbS^{d-1}_{1}}\un_{\{r^{d-1}\tau_{r}(\bfq_{0}, v) \geq \xi\}}F(\bfv)d\bfv - \int^{+\infty}_{\xi}\Phi(x)dx| d\bfq_{0} \\
&=& \int_{\bbT^{d} \cap (\bbR^{d} \setminus \bbQ^{d})}|\int_{\bbS^{d-1}_{1}}\un_{\{r^{d-1}\tau_{r}(\bfq_{0}, v) \geq \xi\}}F(\bfv)d\bfv - \int^{+\infty}_{\xi}\Phi(x)dx| d\bdalp\\
&:=& I_{1} + I_{2}.
\eeqnas
When $\bdalp \in \bigcup^{\infty}_{c \geq \epsilon}D(\kappa, c)$, $\zeta(\bdalp, T) \geq \epsilon^{\frac{1}{\kappa+1}} T^{\frac{1}{\kappa+1}}$, such that 
$$
\zeta(\bdalp,  r^{-\frac{d-1}{2}})^{\bar{\delta}} \geq  \epsilon^{\frac{1}{\kappa+1}\bar{\delta}}r^{-\frac{d-1}{2(\kappa+1)}\bar{\delta}},
$$
and we take $r$, $\epsilon$ satisfying $\epsilon^{\frac{1}{\kappa+1}\bar{\delta}}r^{-\frac{d-1}{2(\kappa+1)}\bar{\delta}} \geq \|M\|$ and apply Theorem~\ref{m.1} with $\lambda' = F$, 
\beqnas
I_{1} &=&  \int_{\bigcup^{\infty}_{c \geq \epsilon}D(\kappa, c)}|\int_{\bbS^{d-1}_{1}}\un_{\{r^{d-1}\tau_{r}(\bfq_{0}, v) \geq \xi\}}F(\bfv)d\bfv - \int^{+\infty}_{\xi}\Phi(x)dx| d\bdalp\\
&=& O(\int_{\bigcup^{\infty}_{c \geq \epsilon}D(\kappa, c)} \zeta(\bdalp, r^{-\frac{d-1}{2}})^{-\bar{\delta}}d\bdalp  + r^{\bar{\kappa}})\\
&=& O(\epsilon^{-\frac{1}{\kappa+1}\bar{\delta}}r^{\frac{d-1}{2(\kappa+1)}\bar{\delta}} + r^{\bar{\kappa}}).
\eeqnas
\beqnas
I_{2} &=& \int_{\bbT^{d} \setminus \bigcup^{\infty}_{c \geq \epsilon}D(\kappa, c)}|\int_{\bbS^{d-1}_{1}}\un_{\{r^{d-1}\tau_{r}(\bfq_{0}, v) \geq \xi\}}F(\bfv)d\bfv - \int^{+\infty}_{\xi}\Phi(x)dx| d\bdalp\\
&\ll& m_{\bbT^{d}}\big(\bbT^{d} \setminus \bigcup^{\infty}_{c \geq \epsilon} D(\kappa, c) \big) \ll \epsilon.
\eeqnas
Taking $\epsilon = r^{\frac{(d-1)\bar{\delta}}{2(\kappa+1 + \bar{\delta})}}$, we derive that for $r^{\frac{(d-1)\bar{\delta}}{2(\kappa+1 + \bar{\delta})}} \leq \|M\|^{-1}$, 
\beqna\label{est.z.2}
\int_{\cF}|\int_{\bbS^{d-1}_{1}}\un_{\{r^{d-1}\tau_{r}(\bfq_{0}, v) \geq \xi\}}F(\bfv)d\bfv - \int^{+\infty}_{\xi}\Phi(x)dx| d\bfq_{0} \ll r^{\frac{(d-1)\bar{\delta}}{2(\kappa+1 + \bar{\delta})}} + r^{\bar{\kappa}}.
\eeqna
Then \eqref{est.z.2} and \eqref{est.fr} lead to
\beqnas
 |\Lambda((\bfQ, \bfV) \in T^{1}(r^{d-1}\cK_{r}): \cT_{r}(\bfQ, \bfV) \geq \xi )- \int^{+\infty}_{\xi}\Phi(x)dx| = O( r^{\frac{(d-1)\bar{\delta}}{2(\bar{\delta} + \kappa +1)}} + r^{\bar{\kappa}} + r^{(d-1)\gamma}).
\eeqnas
Thus we finish the proof of \eqref{est.main.II} by taking $\beta = \min\{\frac{(d-1)\bar{\delta}}{2(\bar{\delta} + \kappa +1)}, \bar{\kappa}, (d-1)\gamma\}$. 

\hfill \qed

%


\end{document}